\documentclass[preprint,12pt]{elsarticle}

\usepackage{amssymb}
\usepackage{amsmath}
\usepackage{algorithm}
\usepackage{algpseudocode}
\usepackage{xcolor}
\usepackage[normalem]{ulem}
\usepackage{caption}
\usepackage{subcaption}
\usepackage{mathtools}

\usepackage{hyperref}
\usepackage{doi}

\usepackage{float}
\newtheorem{remark}{Remark}

\floatstyle{boxed}
\newfloat{customAlg}{h}{ext}
\floatname{customAlg}{Algorithm}

\newcommand{\mtx}[1]{\mathbf{#1}}

\journal{Computers \& Mathematics with Applications}

\begin{document}

\begin{frontmatter}



\title{A Two-Level Direct Solver for the Hierarchical Poincar\'e-Steklov Method}


\author[label1]{Joseph Kump\corref{cor1}} 
\ead{josek97@utexas.edu}
\cortext[cor1]{Corresponding Author}
\author[label1]{Anna Yesypenko} 
\author[label1]{Per-Gunnar Martinsson} 

\affiliation[label1]{organization={Oden Institute, University of Texas at Austin},
            addressline={201 E 24th St}, 
            city={Austin},
            postcode={78712}, 
            state={TX},
            country={USA}}

\begin{abstract}
We introduce a two-level direct solver for the Hierarchical Poincar\'e-Steklov (HPS) method for solving linear elliptic PDEs. 
HPS combines multidomain spectral collocation with a direct solver, enabling high‐order discretizations for highly oscillatory solutions while preserving computational efficiency.
Our method employs batched linear algebra routines with GPU acceleration to reduce the problem to subdomain interfaces, yielding a block-sparse linear system. 
This system is then factorized using a sparse direct solver that employs pivoting to achieve better numerical stability than the original HPS scheme.
For a discretization of local order $p$ involving a total of $N$ degrees of freedom, the initial reduction step has asymptotic complexity $\mathcal O(N p^6)$ in three dimensions.
Nevertheless, the high efficiency of batched GPU routines makes the overall cost for practical purposes independent of polynomial order 
(for order $p=20$ or even higher). 
Additionally, the cost of the sparse direct solver is independent of the polynomial order.
We present a description and justification of our method, along with numerical experiments on three-dimensional problems to evaluate its accuracy and performance.
%
\end{abstract}



\begin{keyword}
sparse direct solvers \sep elliptic PDEs \sep spectral methods \sep domain decomposition \sep GPUs



\end{keyword}

\end{frontmatter}



\section{Introduction}

The manuscript describes a numerical method for accurately and efficiently solving elliptic boundary value problems of the form
\begin{align}
    Au(\mathbf{x}) &= f(\mathbf{x}), &&\mathbf{x} \in \Omega  \nonumber\\
    u(\mathbf{x}) &= g(\mathbf{x}),  &&\mathbf{x} \in \partial \Omega,
    \label{eqn:bvp}
\end{align}
where $\Omega$ is a regular bounded domain in $\mathbb{R}^{d}$ for $d = 2 \text{ or } 3$. Further, $A$ is an elliptic partial differential operator, $g:\partial \Omega \rightarrow \mathbb{R}$ is a given Dirichlet boundary condition, and $f: \Omega \rightarrow \mathbb{R}$ is a given body load. As a concrete example, we focus on the 3D variable-coefficient Helmholtz operator
\begin{equation}
    Au = -\left(\frac{\partial^2}{\partial x_1^2} + \frac{\partial^2}{\partial x_2^2} + \frac{\partial^2}{\partial x_3^2}\right)u - \kappa^2b(\mathbf{x}) u.
    \label{eqn:bvphelmholtz}
\end{equation}
Then (\ref{eqn:bvp}) models time-harmonic wave phenomena, with $\kappa$ denoting the wave-\-number and $b$ a spatially varying, nonnegative function that encodes material heterogeneity.

We introduce a modified solver for the HPS scheme on 2D and 3D domains that forgoes the hierarchical approach of constructing solution operators for subdomain boundaries, instead using a two-level framework. This builds upon work in \cite{yesypenko2022parallel}, which was limited to problems in two dimensions. Rather than recursively constructing a chain of solution operators in an upward pass (as in the original HPS method \cite{2012_spectralcomposite}), we use the DtN maps of each lowest-level subdomain (``leaf box") to assemble a block-sparse system that encodes all leaf box boundaries concurrently. By creating a single larger system, we can interface with multifrontal sparse solvers such as MUMPS, which apply pivoting to improve numerical stability \cite{amestoy2000mumps}. Moreover, the cost of factorizing the sparse system is independent of the choice of subdomain polynomial order $p$. The DtN maps for the subdomains are constructed in parallel with GPU acceleration, where the use of hardware acceleration mitigates the impact of the local polynomial order $p$.

GPU acceleration of PDE solvers has received significant attention in recent years, mostly for low-order numerical methods such as in \cite{dinh2015toward,georgescu2013gpu,klockner2013high, modave2015nodal}. A recurring theme in these works is using GPUs to reduce the problem along each discretization element to a smaller or more structured linear system. We adopt a similar approach, reducing the problem from one on all HPS discretization points to just subdomain boundary points. This yields a significantly smaller sparse system (acting on $O(p^{d-1})$ discretization points per subdomain instead of $O(p^d)$) which is easier to store and faster to factorize, and compatible with black-box direct solvers. The uniformity of this reduction -- applied identically across subdomains -- also enables the use of efficient batched linear algebra routines, allowing us to recompute the required local operators on the fly rather than store them. We assemble the sparse system's data (the DtN maps) and perform localized leaf box solves with batched linear algebra routines on GPUs via PyTorch \cite{imambi2021pytorch}.

Previous work has sought different means to accelerate either the box interior or box boundary solves present in HPS. For instance, \cite{hao2016direct} leverages rank structure in 3D HPS operators to compress matrices and apply fast structured matrix algebra, while \cite{gillman2014n} uses matrix compression to accelerate nested dissection for 2D problems where solutions near the boundary are required. In terms of parallelism, \cite{beams2020parallel} employs OpenMP and MKL to parallelize the upward pass and operator construction within a shared memory framework. More recently, \cite{lucero2024iterative} proposed a distributed memory approach, assembling a global sparse system from all discretization points and solving it iteratively using a GMRES iterative solver with MPI parallelism. In contrast, our method constructs a sparse system using only the box boundary degrees of freedom and retains the use of a direct solver.

This paper provides a summary of our two-level HPS solver. Section~\ref{sec:preliminaries} introduces the block-sparse linear system associated with HPS that we wish to solve, Section~\ref{sec:algorithm} describes the complete algorithm for constructing and solving the reduced system, and Section \ref{sec:numerics} analyzes performance and demonstrates practical applications.

\section{Method Overview}
\label{sec:preliminaries}

In this section we introduce the core components of the method we use to solve the boundary value problem (\ref{eqn:bvp}). We begin by constructing a sparse linear system based on spectral differentiation within subdomains. This leads to a block-sparse global system, which we simplify by eliminating interior degrees of freedom to expose a smaller system that involves only degrees of freedom associated with discretization nodes on subdomain boundaries.

\subsection{Global Discretization}
Equation~(\ref{eqn:bvp}) can be solved using the Hierarchical Poincar\'e-Steklov (HPS) method \cite{2012_spectralcomposite,2013_martinsson_ItI,martinsson2015hierarchical}, which combines a multidomain spectral collocation discretization with a nested dissection based direct solver. The discretization starts with a flat tesselation into subdomains (``boxes"), where the PDE is enforced locally through spectral collocation in each subdomain, and then continuity of the normal derivative is enforced across boundaries. Like other multidomain spectral methods, HPS enables use of high-order local discretizations that mitigate the pollution effect \cite{babuska1995generalized, babuska1997pollution, galkowski2023does}. It also requires only minimal subdomain overlap, unlike high-order finite difference methods, which reduces computational cost in direct solvers \cite[Chapter~24]{martinsson2019fast} and inter-process communication in parallel settings \cite{davis2006direct}. In previous HPS solvers such as \cite{2012_spectralcomposite,hao2016direct,beams2020parallel, chen2024fast}, the resulting discretized system is then solved via the hierarchical merging of Dirichlet-to-Neumann (DtN) maps.
\label{ssec:disc}

\begin{figure}
    \centering
    \includegraphics[width=0.9\linewidth]{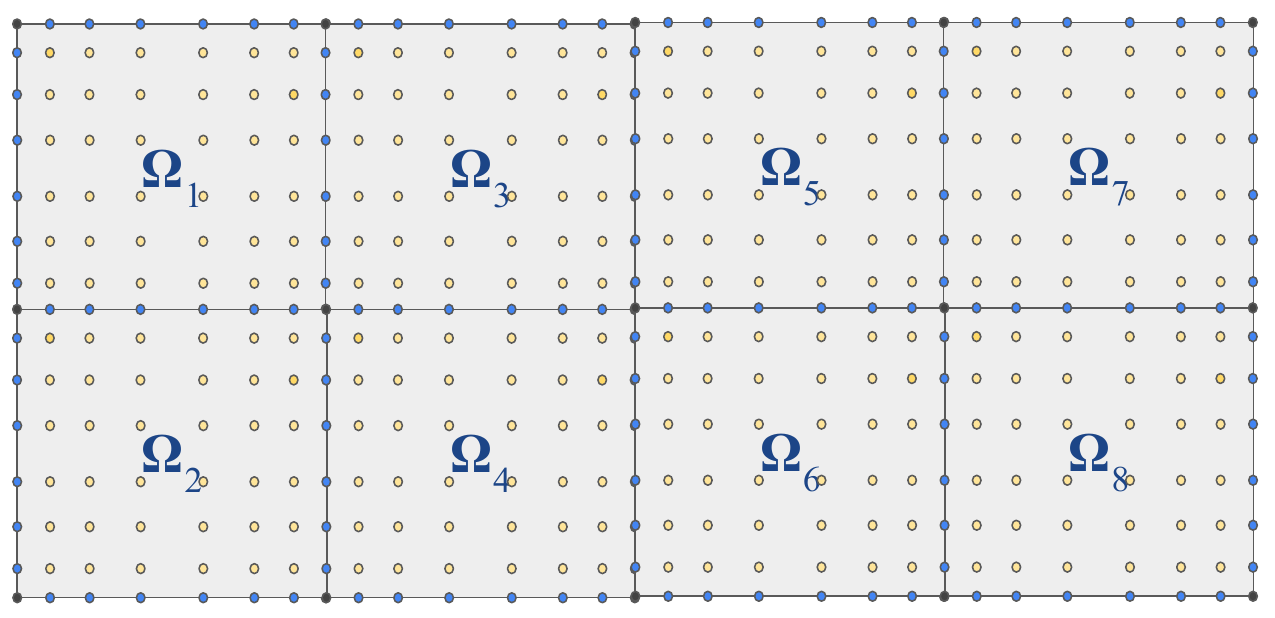}
    \caption{Diagram of a 2D domain partitioned into subdomains for HPS. Yellow nodes are subdomain interiors, blue nodes are subdomain boundaries, and gray nodes are corners.}
    \label{fig:2d-domain}
\end{figure}

The domain $\Omega \subset \mathbb{R}^d$ is partitioned into a series of non-overlapping subdomains, often called ``leaf boxes" or just ``boxes". Generally these subdomains are rectangular and of the same size, but their shape may be modified through the use of parameter maps (as seen in Section~\ref{ssec:curved}). An example of a 2D discretization with uniform rectangular subdomains is shown in Figure~\ref{fig:2d-domain}. On each subdomain $\Omega^{(j)}$ we place a tensor product of $p^d$ Chebyshev nodes, with the choice of polynomial order $p$ consistent across all subdomains. These nodes are shared along the boundaries between subdomains: edges in 2D and faces in 3D. We aim to solve for the solution $u$ of our boundary value problem on these discretization points, which we denote as the vector $\mathbf{u}$.

On the subdomain interiors, we locally enforce the differential equation (\ref{eqn:bvp}) using a spectral collocation method \cite{2000_trefethen_spectral_matlab}. On the subdomain boundaries, we enforce continuity of the Neumann derivative of the two adjacent subdomains. More detailed formulations for how we encode these conditions can be found in \cite[Chapter~24]{martinsson2019fast}. Together, the matrices for interior and boundary solves form the block rows of a larger system $\mathbf{A}$ that concurrently solves for all subdomain interior and boundary points of our discretization, 
\begin{equation}
\label{eq:Au=f}
\mathbf{A} \mathbf{u} = \mathbf{f}.
\end{equation}
For a discretization node $\mathbf{x}_i \in \Omega$, the value of $[\mathbf{Au}](i)$ is
\begin{equation}
\label{eq:Aentries}
    [\mathbf{Au}](i) \approx \begin{cases}
        Au(\mathbf{x}_i) \text{ if } \mathbf{x}_i \in \text{ interior of some box } \Omega^{(j)}, \\
        \left. \frac{\partial u}{\partial n}(\mathbf{x}_i)\right|_{\Omega^{(j)}} + \left. \frac{\partial u}{\partial n}(\mathbf{x}_i)\right|_{\Omega^{(k)}} \text{ if } \mathbf{x}_i \in \text{ boundary } \partial \Omega^{(j)} \cap \partial \Omega^{(k)}. \\
    \end{cases}
\end{equation}
Nodes that lie on corners (and edges in 3D) in principle require special treatment, but in most applications they can surprisingly be dropped from consideration entirely. See Remark \ref{remark:corner}.

\begin{figure}
\centering
Matrix Sparsity Patterns
\begin{subfigure}{.5\textwidth}
  \centering
  \includegraphics[width=\textwidth,trim={0cm 0cm 0cm 0cm},clip]{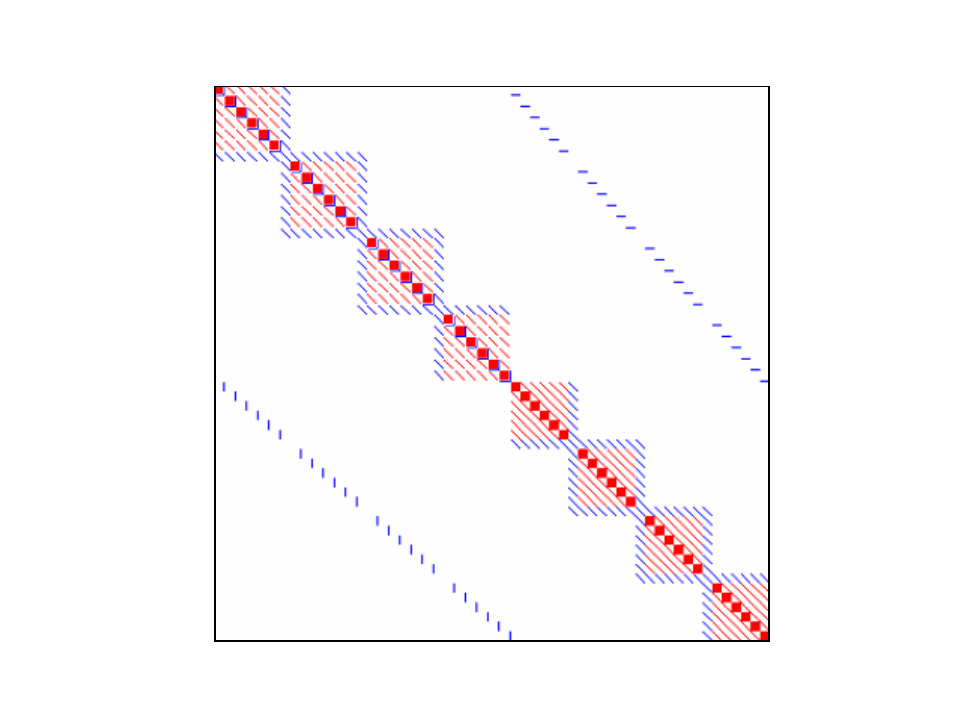}
  \caption{Standard blockwise ordering.}
  \label{fig:a-naive-denoted}
\end{subfigure}%
\begin{subfigure}{.5\textwidth}
  \centering
  \includegraphics[width=\textwidth,trim={0cm 0cm 0cm 0cm},clip]{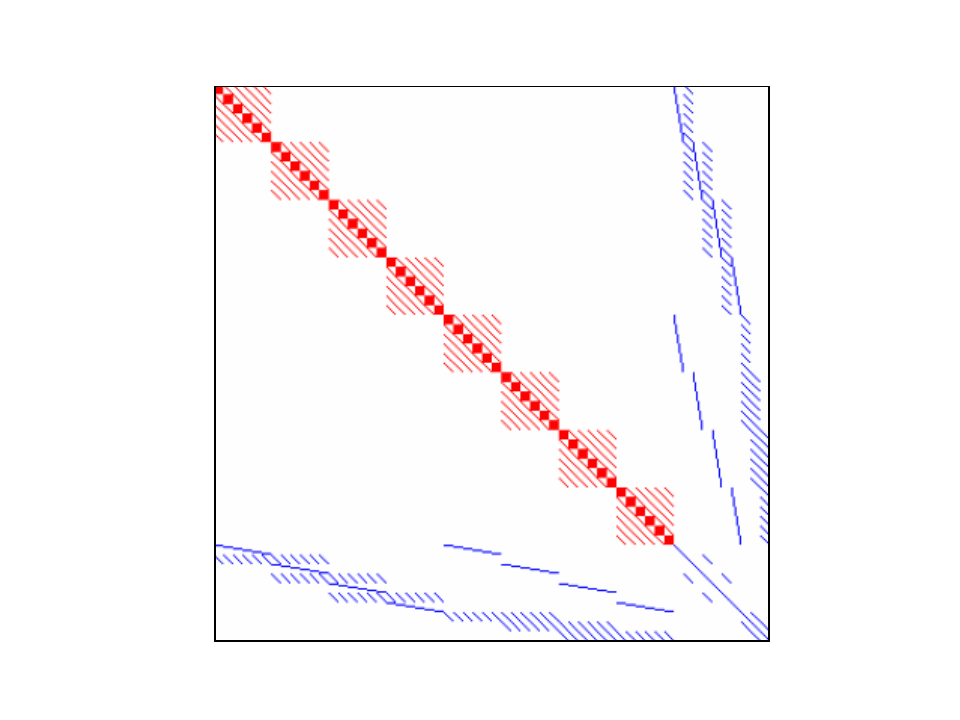}
  \caption{Permuted so that interior nodes go first.}
  \label{fig:a-rearranged-denoted}
\end{subfigure}
\caption{Sparsity pattern of $\mathbf{A}$ for a 2D domain of $4 \times 2$ subdomains as in Figure \ref{fig:2d-domain}. Entries associated with interior nodes are red, and entries associated with boundary nodes are blue. Section \ref{ssec:static-condensation} provides details.}
\label{fig:a-naive}
\end{figure}

The matrix $\mathbf{A}$ is block-sparse with distinctive structures that are illustrated in Figure \ref{fig:a-naive}. It contains relatively dense blocks corresponding to the interactions within subdomain interiors (red in Fig. \ref{fig:a-naive-denoted}) bordered by smaller submatrices representing the Neumann derivative approximations between ``left" and ``right" adjacent subdomains (blue in Fig. \ref{fig:a-naive-denoted}). The interactions between ``up" and ``down'' adjacent subdomains are represented by the blue submatrices far from the main diagonal.

\begin{remark}[Corner and edge nodes]
\label{remark:corner}
The formula (\ref{eq:Aentries}) does not specify how to handle nodes that lie on corners and edges. Conveniently, it turns out that these nodes can be simply dropped from consideration in any situation where the differential operator in  (\ref{eqn:bvp}) does not involve cross terms of the form $\partial_{i}\partial_{j}u$ for $i \neq j$. In situations where cross terms do appear, one reinterpolates the $p\times p$ Chebyshev grid on each face to a $(p-1)\times (p-1)$ grid of Legendre nodes. This way the global system again will involve no discretization nodes on edges or corners. For details, see \cite{hao2016direct,martinsson2015hierarchical}.
\end{remark}

\subsection{Static condensation: Elimination of interior nodes}
\label{ssec:static-condensation}

\begin{figure}[t]
\centering
\includegraphics[width=350pt]{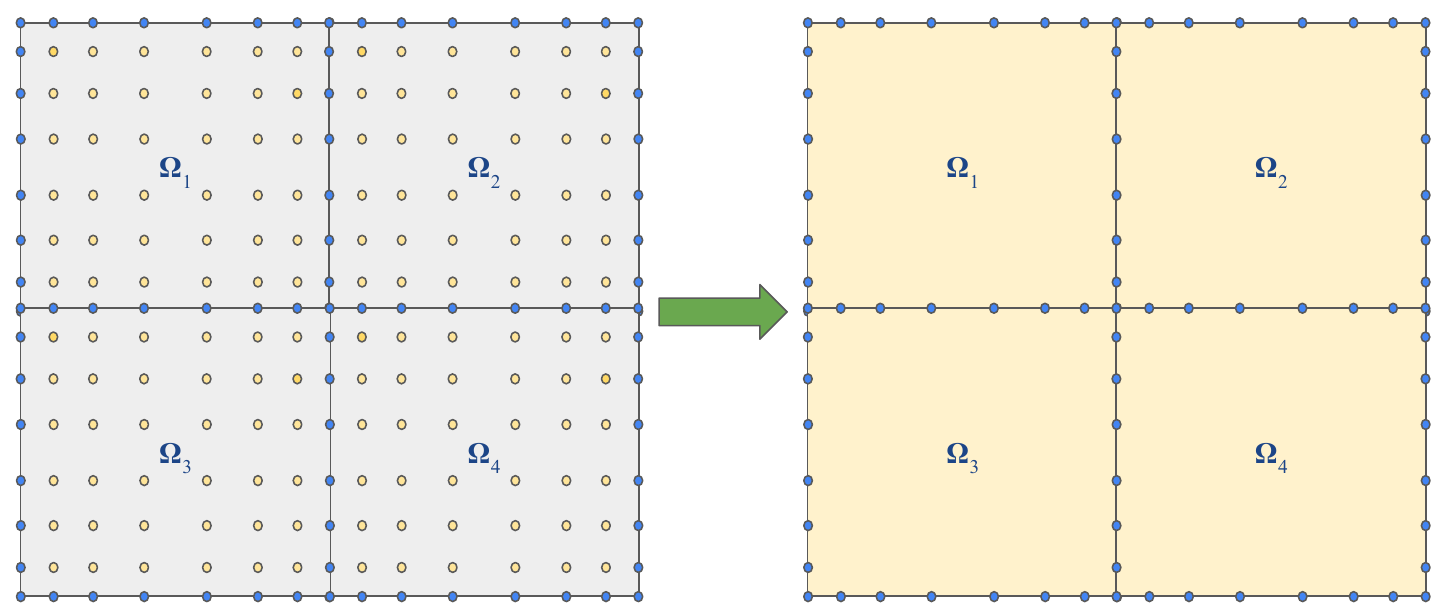}
\caption{We reduce the problem from one on the entire discretized domain to one on subdomain boundaries only.}
\label{fig:staticCondensation}
\end{figure}

For computational efficiency, we in practice never explicitly form the large sparse matrix $\mtx{A}$. Instead we perform a pre-computation on each box that eliminates all nodes that are interior to the box, to form a sparse matrix $\mtx{T}$ that encodes only interactions between boundary nodes illustrated on a domain in Figure~\ref{fig:staticCondensation}. The matrix $\mtx{T}$ is much smaller than $\mtx{A}$ (roughly by a factor of $p$), but is less sparse.

The first step is to partition the discretization nodes into two sets: one that holds all nodes that are interior to a box, and one that holds all nodes that sit on domain boundaries. Let $I_{\rm i}$ and $I_{\rm b}$ be two index vectors that identify the two sets. (In other words, $I_{\rm i}$ marks all diagonal entries that are red in Figure \ref{fig:a-naive}, while $I_{\rm b}$ marks all blue diagonal entries.) Then we partition the linear system (\ref{eq:Au=f}) accordingly, to obtain
\begin{equation}
\label{eq:blocksys1}
\left[\begin{array}{rr}
\mtx{A}_{\rm ii} & \mtx{A}_{\rm ib} \\
\mtx{A}_{\rm bi} & \mtx{A}_{\rm bb} 
\end{array}\right]
\left[\begin{array}{rr}
\mtx{u}_{\rm i} \\ \mtx{u}_{\rm b}
\end{array}\right]
=
\left[\begin{array}{rr}
\mtx{f}_{\rm i} \\ \mtx{f}_{\rm b}
\end{array}\right],
\end{equation}
where
$$
\mtx{A}_{\rm ii} = \mtx{A}(I_{\rm i},I_{\rm i}),\qquad
\mtx{A}_{\rm ib} = \mtx{A}(I_{\rm i},I_{\rm b}),\qquad
\mtx{A}_{\rm bi} = \mtx{A}(I_{\rm b},I_{\rm i}),\qquad
\mtx{A}_{\rm bb} = \mtx{A}(I_{\rm b},I_{\rm b}),
$$
and where
$$
\mtx{u}_{\rm i} = \mtx{u}(I_{\rm i}),\qquad
\mtx{u}_{\rm b} = \mtx{u}(I_{\rm b}),\qquad
\mtx{f}_{\rm i} = \mtx{f}(I_{\rm i}),\qquad
\mtx{f}_{\rm b} = \mtx{f}(I_{\rm b}).
$$
A key observation here is that the matrix $\mtx{A}_{\rm ii}$ is block diagonal, since the interior nodes of any single box do not directly communicate with the interior nodes of any other box. This means that we can apply $\mtx{A}_{\rm ii}^{-1}$ to vectors via local computations on each box that are completely unconnected. Another key observation is that even for variable coefficient differential operators, the matrices $\mtx{A}_{\rm bi}$ and $\mtx{A}_{\rm bb}$ both consist of repeated, identical submatrices, at least in the case where all subdomains have uniform $p$. This is because the differential operators they numerically represent - specifically the Neumann derivative to ensure Neumann continuity - are universal across all boxes in the domain. We exploit this invariance to accelerate many of the box computations. A sparsity plot of $\mathbf{A}$ partitioned as in (\ref{eq:blocksys1}) is shown in Figure \ref{fig:a-rearranged-denoted}.

\begin{figure}
Factorized Matrix Sparsity Patterns
\centering
\begin{subfigure}{.24\textwidth}
  \centering
  $\tilde{\mathbf{A}}$
  \includegraphics[width=\textwidth,trim={2.5cm 0cm 2.5cm 0cm},clip]{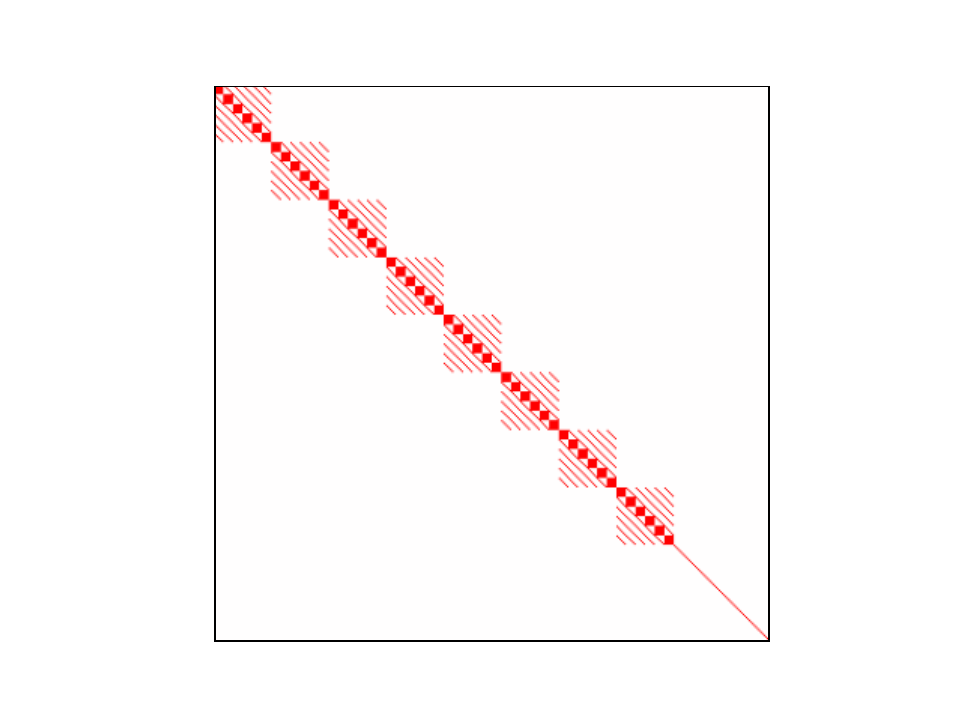}
  \caption{$\mathbf{A}_{\rm ii}$ in red.}
  \label{fig:aii-values}
\end{subfigure}%
\begin{subfigure}{.24\textwidth}
  \centering
  $\mathbf{L}$
  \includegraphics[width=\textwidth,trim={2.5cm 0cm 2.5cm 0cm},clip]{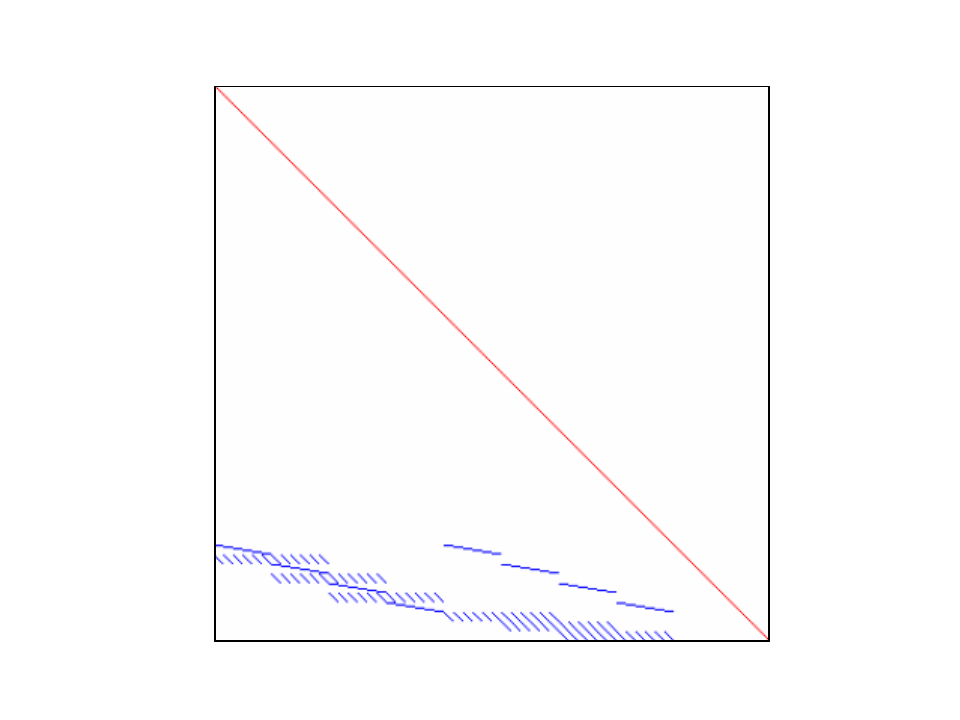}
  \caption{$\mtx{A}_{\rm bi}$ in blue.}
  \label{fig:l-values}
\end{subfigure}%
\begin{subfigure}{.24\textwidth}
  \centering
  $\mathbf{D}$
  \includegraphics[width=\textwidth,trim={2.5cm 0cm 2.5cm 0cm},clip]{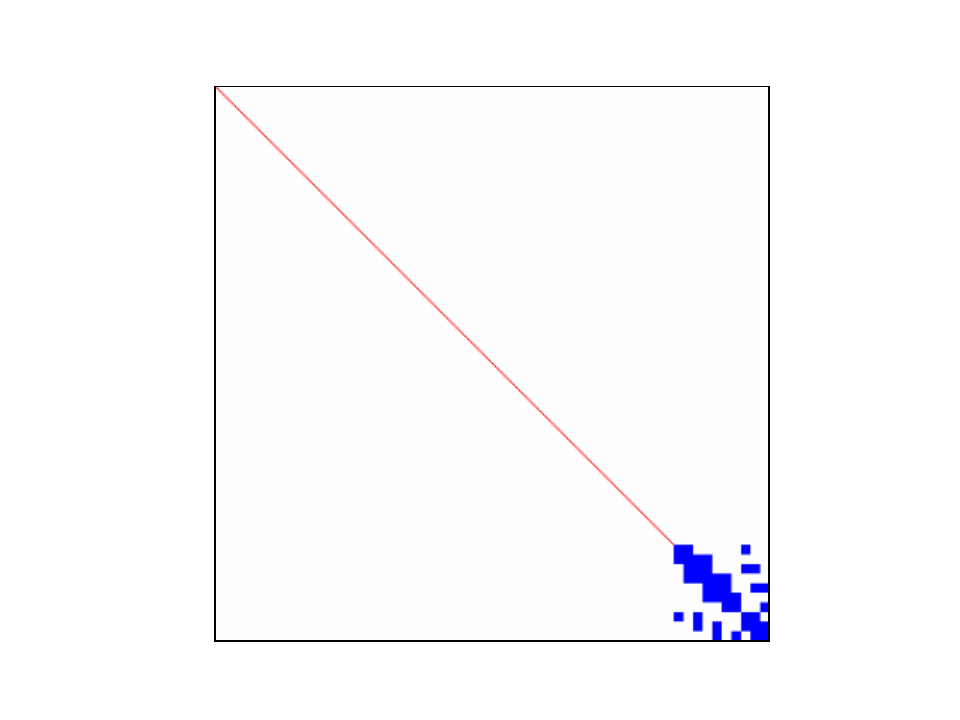}
  \caption{$\mathbf{T}$ in blue.}
  \label{fig:d-values}
\end{subfigure}
\begin{subfigure}{.24\textwidth}
  \centering
  $\mathbf{U}$
  \includegraphics[width=\textwidth,trim={2.5cm 0cm 2.5cm 0cm},clip]{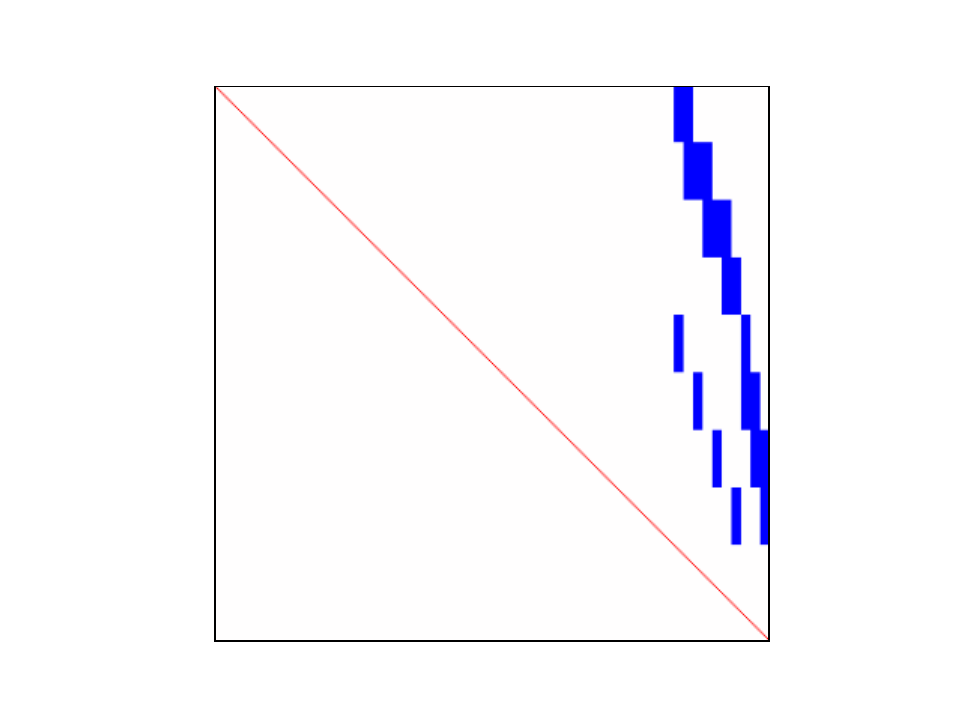}
  \caption{$-\mathbf{S}$ in blue.}
  \label{fig:u-values}
\end{subfigure}
\caption{Sparsity patterns of (with $\tilde{\mathbf{A}}$)$\mathbf{LDU}$ factorization for permuted $\mathbf{A}$. This is for a 2D problem with $4 \times 2$ subdomains as in Figure~\ref{fig:a-naive}. Blocks corresponding to particular submatrices in the factorization are highlighted red or blue.}
\label{fig:ldu}
\end{figure}

The solver we use relies on the following factorization of the coefficient matrix in (\ref{eq:blocksys1}):
\begin{equation}
\label{eq:factor1}
\mtx{A} = 
\underbrace{\left[\begin{array}{rr}
\mtx{A}_{\rm ii} & \mtx{0}\\
\mtx{0} & \mtx{I}
\end{array}\right]}_{=\tilde{\mtx{A}}}
\underbrace{\left[\begin{array}{rr}
\mtx{I} & \mtx{0}\\
\mtx{A}_{\rm bi} & \mtx{I}
\end{array}\right]}_{=\mtx{L}}
\underbrace{\left[\begin{array}{rr}
\mtx{I} & \mtx{0}\\
\mtx{0} & \mtx{A}_{\rm bb} - \mtx{A}_{\rm bi}\mtx{A}_{\rm ii}^{-1}\mtx{A}_{\rm ib}
\end{array}\right]}_{=\mtx{D}}
\underbrace{\left[\begin{array}{rr}
\mtx{I} & \mtx{A}_{\rm ii}^{-1}\mtx{A}_{\rm ib}\\
\mtx{0} & \mtx{I}
\end{array}\right]}_{=\mtx{U}}
\end{equation}
That (\ref{eq:factor1}) holds is easily established by simply multiplying the factors together. It will be convenient to introduce two new matrices
$$
\mtx{S} := -\mtx{A}_{\rm ii}^{-1}\mtx{A}_{\rm ib}
\qquad\mbox{and}\qquad
\mtx{T} := 
\mtx{A}_{\rm bb} - \mtx{A}_{\rm bi}\mtx{A}_{\rm ii}^{-1}\mtx{A}_{\rm ib} = 
\mtx{A}_{\rm bb} + \mtx{A}_{\rm bi}\mtx{S}.
$$
From (\ref{eq:factor1}), we then easily obtain the inversion formula
\begin{equation}
\label{eq:factor2}
\mtx{A}^{-1} = 
\underbrace{\left[\begin{array}{rr}
\mtx{I} & \mtx{S}\\
\mtx{0} & \mtx{I}
\end{array}\right]}_{=\mathbf{U}^{-1}}
\underbrace{\left[\begin{array}{rr}
\mtx{I} & \mtx{0}\\
\mtx{0} & \mtx{T}^{-1}
\end{array}\right]}_{=\mathbf{D}^{-1}}
\underbrace{\left[\begin{array}{rr}
\mtx{I} & \mtx{0}\\
-\mtx{A}_{\rm bi} & \mtx{I}
\end{array}\right]}_{=\mathbf{L}^{-1}}
\underbrace{\left[\begin{array}{rr}
\mtx{A}_{\rm ii}^{-1} & \mtx{0}\\
\mtx{0} & \mtx{I}
\end{array}\right]}_{=\tilde{\mathbf{A}}^{-1}}.
\end{equation}
Observe that the matrix $\mtx{T}$ is sparse and is much smaller than the original matrix $\mtx{A}$. However, it has a significantly larger ratio of non-zero elements per row than $\mtx{A}$. The submatrices $\mtx{S}$, $\mtx{A}_{\rm bi}$, and $\mathbf{A}_{\rm ii}^{-1}$ are all block matrices, with $\mathbf{A}_{\rm bi}$ in particular consisting of the same repeated block.

To compute the solution $\mtx{u}$ to the linear system (\ref{eq:Au=f}) we now need to execute four steps, each one corresponding to one factor in (\ref{eq:factor2}):

\begin{enumerate}
\item Form the vector $\mtx{v}_{\rm i} = \mtx{A}_{\rm ii}^{-1}\mtx{u}_{\rm i}$.
Exploit that $\mtx{A}_{\rm ii}$ is block diagonal to do this through independent local computations, executed via batched linear algebra.
\item Form the vector $\mtx{g}_{\rm b} = \mtx{f}_{\rm b} - \mtx{A}_{\rm bi}\mtx{v}_{\rm i}$. We observe that $\mtx{A}_{\rm bi}$ is very sparse and specifically consists of one repeated block, so the vector $\mtx{g}_{\rm b}$ can again be formed via a repeated local computation on each box that is accelerated via GPUs (this is a ``matrix free'' approach).
\item Assemble the block sparse matrix $\mtx{T}$ through a set of local computations on all boxes. Then solve the global linear system $\mtx{T}\mtx{u}_{\rm b} = \mtx{g}_{\rm b}$ using a highly optimized sparse direct solver.
\item Once $\mtx{u}_{\rm b}$ is available, we can reconstruct the full solution via the formula $\mtx{u}_{\rm i} = \mtx{v}_{\rm i} + \mtx{S}\mtx{u}_{\rm b}$. 
\end{enumerate}

In situations where the system (\ref{eqn:bvp}) needs to be solved for multiple right-hand sides, further acceleration can be obtained by pre-computing and storing an LU factorization of the matrix $\mtx{T}$. Details of how the algorithm was implemented efficiently are provided in Section \ref{sec:algorithm}. 

\subsection{Discussion}

In this manuscript we described the solver from a linear algebraic point of view, as this helps clarify how we attain acceleration via batched linear algebra. However, all steps we take have direct analogues to steps in the original HPS method, as described in \cite[Chapter~25]{martinsson2019fast}. To be precise, the blocks of $\mtx{T}$ contain the Dirichlet-to-Neumann operators that are formed in HPS. The matrix $\mtx{S}$ holds the ``solution operators'' that reconstruct the solution to (\ref{eqn:bvp}) at the interior nodes, once the solution at the boundaries has been identified in the global sparse solve.

In the original HPS method the sparse matrix $\mtx{T}$ is never formed explicitly. Instead, the same solution process is continued hierarchically: first a set of coarser boxes is created by taking the lowest-level boxes and merging them into pairs. This creates a new partition between ``interior'' and ``boundary'' boxes which is used to further factorize the matrix $\mtx{T}$ itself to form an even smaller block sparse matrix $\mtx{T}'$, and so on. The solution maps that make up these levels of matrices ($\mtx{T}$, $\mtx{T}'$, and ongoing) can be formed in parallel within each level.

Highly-optimized multifrontal solvers like MUMPS can use pivoting techniques such as partial threshold pivoting, relaxed pivoting \cite{duff2007towards}, and tournament pivoting \cite{grigori2008communication} to improve numerical stability of a matrix factorization like sparse LU \cite{amestoy2024mumps}. In theory such pivoting techniques could be utilized on matrices present in previous approaches to HPS to improve their stability. However, encoding solutions for all shared faces into one sparse matrix operator $\mathbf{T}$ allows such pivoting techniques to be applied holistically to a solve for $\mathbf{u}$ on all such points. This means a wider range of pivots can be applied, improving the factorization's stability while preserving an efficient sparsity pattern. In effect we ``look at" a larger portion of the problem at once, allowing us (through the use of a multilevel solver) to perform more optimizations in the factorization process. This approach also enables us to use highly optimized standard libraries in the construction and factorization of $\mathbf{T}$.

In addition, the non-hierarchical nature of our assembly of $\mathbf{T}$ -- relying only on the lowest-level DtN maps and not higher level DtN or solution operators -- gives this approach uniformity across the computation of its blocks which is helpful for parallelization. Individually the DtN maps that make up $\mtx{T}$ are relatively dense matrices, but their entries can be computed on GPUs if available via batched linear algebra routines \cite{dongarra2017optimized}. This allows us to use GPU accelerators on a significant portion of our method's computations.

\section{Algorithm Description}
\label{sec:algorithm}

In this section we elaborate on the  implementation of our HPS solver, highlighted by the four steps presented in Section \ref{ssec:static-condensation}. We discuss the numerical challenges present in both the local box operations and the formulation of a system for the reduced problem on domain boundaries.

\subsection{Local Computations}

As shown in Section \ref{ssec:static-condensation}, our HPS solver consists of four steps corresponding to the factors $\tilde{\mtx{A}}^{-1}$, $\mtx{L}^{-1}$, $\mtx{D}^{-1}$, and $\mtx{U}^{-1}$. We can avoid forming these matrices explicitly, and instead apply them in our solvers through matrix-free approaches. The steps to accomplish this fit into two categories of computations.

The first set are batched linear algebra routines that represent $\tilde{\mtx{A}}, \mtx{L},$ and $\mtx{U}$: the factors are primarily made of submatrices that are block-sparse in structure, such as $\mtx{A}_{\rm ii}^{-1}$, $\mtx{S}$, and $\mtx{A}_{\rm bi}$. Each of these matrices have a limited number of nonzero blocks per row (with each block-row corresponding to solving for $\mtx{u}$ on one box $\Omega^{\tau}$) that are only in specific indices relative to the row index. Thus we stack these nonzero blocks into three-dimensional arrays of size $(\text{number of blocks } \times \text{ size of each block})$, with their indices within the matrix inferred from the order in which they are stored. We then apply them to our vector data in a matrix-free approach (using ``tensor" objects and operations in PyTorch). Practically speaking, we are representing the action of the matrix-vector product with a series of smaller matrix-vector products corresponding to each nonzero block. By doing this we can use batched linear algebra routines (also present in PyTorch) for the matrix-vector multiplications, which reduce the memory cost and computational time.

Further optimizations can be applied to specific submatrices. For instance, since $\mtx{A}_{\rm ii}$ is block-diagonal, we can apply $\mtx{A}_{\rm ii}^{-1}$ to a vector by simply performing a number of independent small system solves. Meanwhile, since $\mtx{A}_{\rm bi}$ represents part of a Neumann operator used to enforce Neumann continuity on boundary nodes, it is unaffected by our differential operator $A$ so its blocks $\mtx{A}_{\rm bi}^{\tau}$ are identical across all boxes ${\Omega^{\tau}}$ of our discretization. This means $\mtx{A}_{\rm bi}$ does not require a stacked 3-dimensional array to represent its multiplication; we may represent it as a single small block applied to each of the subvectors of $\mtx{u}$ representing a box. This reduces the memory footprint even further than computations with stacked 3D arrays.

The exception to this batched approach is the sparse matrix $\mtx{T}$ that solves for $\mtx{u}$ on boundary nodes. While $\mtx{T}$ is also block-sparse, it must be inverted (or, more accurately, factorized for use in a solver) to apply $\mtx{D}^{-1}$, and unlike $\mtx{A}_{\rm ii}$ it is not block-diagonal. Thus there is no similarly convenient means to invert its components independently, and we instead interface with MUMPS to produce a factorization of it for use in our method.

However, while $\mtx{T}^{-1}$ cannot easily be applied with a batched matrix-free method in a comparable way to the other submatrices, we may assemble $\mtx{T}$ in batch before factorizing it. We can write $\mtx T$ as
\begin{equation}
    \label{eqn:t-block}
    \mtx{T} = \mtx{A}_{\rm bb} + \mtx{A}_{\rm bi} \mtx{S} = \begin{bmatrix}
    \mtx{A}_{\rm bb} & \mtx{A}_{\rm bi}
\end{bmatrix} \begin{bmatrix}
    \mtx{I}\\
    \mtx{S}
\end{bmatrix}.
\end{equation}
Since $\mtx{S}$ itself consists of blocks corresponding to each box $\Omega^{\tau}$ that can be neatly represented by a 3D array, we may assemble and store the blocks of $\mtx{T}$ in a similar fashion, and then apply the blocks of $\left[\mtx{A}_{\rm bb} | \mtx{A}_{\rm bi}\right]$ to this array (like $\mtx{A}_{\rm bi}$, $\mtx{A}_{\rm bb}$ is identical across blocks so we can represent it as a single small matrix). This gives us the blocks $\mtx{T}^{\tau}$ that make up $\mtx{T}$. From here we can then construct row and column index arrays that align each $\mtx{T}^{\tau}$ to the subsets of $\mtx{u}_b$ that reside on $\partial \Omega^{\tau}$, then construct $\mtx{T}$ in a sparse format and factorize it with MUMPS, in our case accessed via petsc4py.

\begin{remark}[Relation to previous HPS implementations] 
    The submatrices present in our $\tilde{\mtx{A}} \mtx{L} \mtx{D} \mtx{U}$ factorization of $\mtx{A}$ each represent operators present in previous formulations of HPS. In particular, $\mtx{A}_{\rm bi}$ and $\mtx{A}_{\rm bb}$ are differential matrices approximating the Neumann derivative, the blocks of $\mtx{S}$ are the solution operators for solving box interiors, and the blocks of $\mtx{T}$ are the DtN maps consisting of composing Neumann differentiation on solution operators to enforce Neumann continuity. Readers previously familiar with HPS likely notice (\ref{eqn:t-block}) resembles the formula for constructing a DtN map.
\end{remark}

\subsection{Build and Solve stages}


Similar to previous HPS implementations, our method can be divided into a build stage (assembling and factorizing the operator $\mtx{T}$) and a solve stage (using the LU factors of $\mtx{T}$, along with other operators, to solve for $\mtx{u}$).

\begin{customAlg}
    \begin{minipage}{135mm}
    \begin{center}
    \textsc{Build stage}
    \end{center}
    \begin{enumerate}
    \item Initialize a sparse matrix $\mtx{T} = \mtx{0}$.
    \item \textit{For each box $\tau$, compute the DtN matrix $\mtx{T}^{\tau}$ as described in Section 3.1, and add it to the relevant slot of $\mtx{T}$.}\\
    \textbf{for} (all boxes $\tau$) \hfill\fbox{\textit{In parallel.}}
    \begin{itemize}
    \item Form the local equilibrium matrix $\mtx{A}^{\tau}$.
    \item Form the local solution operator $\mtx{S}^{\tau} = -(\mtx{A}_{\rm ii}^{\tau})^{-1}\mtx{A}_{\rm ib}^{\tau}$
    \item Form the local DtN matrix $\mtx{T}^{\tau} = \mtx{A}_{\rm bi}^{\tau}\mtx{S}^{\tau}$.
    \item Add the matrix $\mtx{T}^{\tau}$ to the relevant block of $\mtx{T}$.
    \end{itemize}
    \textbf{end for}
    \item Form the $\mtx{LU}$ factorization of $\mtx{T}$.
    \end{enumerate}
    \end{minipage}
    \caption{constructs and factorizes $\mtx{T}$.}
    \label{alg:buildStage}
\end{customAlg}

An outline of the build stage is shown in Algorithm~\ref{alg:buildStage}. In short, we construct the DtN map for each subdomain $\tau$ in our partitioned domain, then assemble $\mtx{T}$ out of these DtNs. Each DtN is of the same size ($\approx 4p \times 4p$ in 2D) and is computed using the same formulation of differentiation operators. Thus they can be stored in ``stacked" three-dimensional arrays and effectively assembled using batched linear algebra routines such as those in PyTorch.

However, a challenge with utilizing GPUs for DtN assembly is memory storage. In 3D each DtN map is size $\approx 6p^2 \times 6p^2$. Building them requires solving matrix systems of the spectral differentiation operators on interior points, which are size $\approx (p-2)^3 \times (p-2)^3$. These matrices grow quickly as $p$ increases (at $p=20$ they are $\approx 5.7$ million and $34$ million entries, respectively). For larger problem sizes it is impossible to store all DtN maps on a GPU concurrently. The interior differentiation submatrices in particular do not have to be stored beyond formulating their box's DtN map, but they are much larger for high $p$ and present considerable memory overhead. To handle this storage requirement we have implemented a two-level scheduler that stores a limited number of DtNs on the GPU at once, and assembles a smaller number of DtNs in batch at a time. For $p \geq 16$ on a V100 GPU the best results occur when only one DtN is assembled at a time. However, we may still store multiple DtNs on the GPU before transferring them collectively to the CPU.

Once the DtNs $\{\mtx{T}^{\tau}\}$ are assembled, we add them to the relevant blocks of $\mtx{T}$, which is represented as a sparse matrix on the CPU. The row and column indexing arrays used to determine these blocks are also constructed in parallel with PyTorch. Then we factorize $\mtx{T}$ using MUMPS. This factorization can be reused for multiple problems with the same differential operator $A$, avoiding the need to rebuild it.

\begin{customAlg}
    \begin{minipage}{135mm}
    \begin{center}
    \textsc{Solve stage}
    \end{center}
    \begin{enumerate}
    \item \textit{Initialize the solution vector:}
    $\mtx{u} = \mtx{0}$.
    \item \textit{In a loop over all boxes, form the local solutions $\mtx{v}^{\tau}$ and the local equivalent load vectors $\mtx{g}^{\tau}$:}\\
    \textbf{for} (all boxes $\tau$) \hfill\fbox{\textit{In parallel.}}
    \begin{itemize}
    \item Form the matrices $\mtx{A}_{\rm ii}^{\tau}$ and $\mtx{A}_{\rm bi}^{\tau}$
    \item Form $\mtx{u}(I_{\rm i}^{\tau}) = (\mtx{A}^{\tau}_{\rm ii})^{-1}\mtx{f}(I^{\tau}_{\rm i})$.
    \hfill \textit{(So $\mtx{u}(I_{\rm i}^{\tau})$ holds $\mtx{v}^{\tau}$.)}
    \item Form $\mtx{f}(I^{\tau}_{\rm b}) = \mtx{f}(I^{\tau}_{\rm b}) - \mtx{A}^{\tau}_{\rm bi}\mtx{u}(I_{\rm i}^{\tau})$.
    \hfill \textit{(So $\mtx{f}(I_{\rm b}^{\tau})$ holds $\mtx{g}^{\tau}$.)}
    \end{itemize}
    \textbf{end for} 
    \item \textit{Solve the global system using a sparse direct solver with the precomputed $\mtx{LU}$ factorization of $\mtx{T}$:} 
    $\mtx{u}(I_{\rm b}) = \mtx{U}^{-1}\mtx{L}^{-1}\mtx{f}(I_{\rm b})$ 
    \item \textit{Fill in the solution at all the interior nodes:}\\
    \textbf{for} (all boxes $\tau$) \hfill\fbox{\textit{In parallel.}}
    \begin{itemize}
    \item 
    $\mtx{u}(I_{\rm i}^{\tau}) = \mtx{u}(I_{\rm i}^{\tau}) + \mtx{S}^{\tau}\mtx{u}(I_{\rm b}^{\tau})$.
    \end{itemize}
    \textbf{end for} 
    \end{enumerate}
    \end{minipage}
    \caption{solve for the PDE given a factorized system $\mtx{T}$.}
    \label{alg:solveStage}
\end{customAlg}

Algorithm~\ref{alg:solveStage} outlines the solve stage. Here we compute $\mtx{u}$ given a differential operator $A$, resulting boundary system $\mtx{T}$, Dirichlet data, and a body load $\mtx{f}$. We must reduce $\mtx{f}$ to the subdomain boundaries before applying it to solve for $\mtx{u}$ on box boundary points as described in Section \ref{ssec:static-condensation}. In addition, after solving for boundary $\mtx{u}$ we must apply local solution operators $\mtx{S}^{\tau}$ to $\mtx{u}$ on each subdomain boundary to obtain $\mtx{u}$ on the subdomain interiors. Each of these operations can be done in parallel using blocks of $\mtx{A}_{bi}, \mtx{A}_{ii}^{-1}$, and $\mtx{S}$. Although these matrices are needed to assemble the local DtN maps, we do not store them after the build stage. For large problems, GPU memory is insufficient, and CPU-GPU transfers are costly, thus it is more efficient to recompute the required block operators during each solve.

The uniformity of DtN assembly and interior subdomain solves is clear when one choice of $p$ is used across all subdomains. However, we could instead select $p$ from a small fixed set of values (perhaps $\approx 6$ different choices) to resemble adaptive methods and still retain much of the benefit. In this case we would have separate arrays storing the DtNs for each choice of $p$, but as long as enough subdomains use each $p$ the corresponding submatrix multiplications can still be batched. The same infrastructure that enables Chebyshev-Legendre interpolation (see Remark \ref{remark:corner}) can also enable interpolation between adjacent boxes with different $p$ in the assembly of $\mtx{T}$, which will be finalized in future work.

\begin{figure}[t]
\centering
\includegraphics[width=250pt]{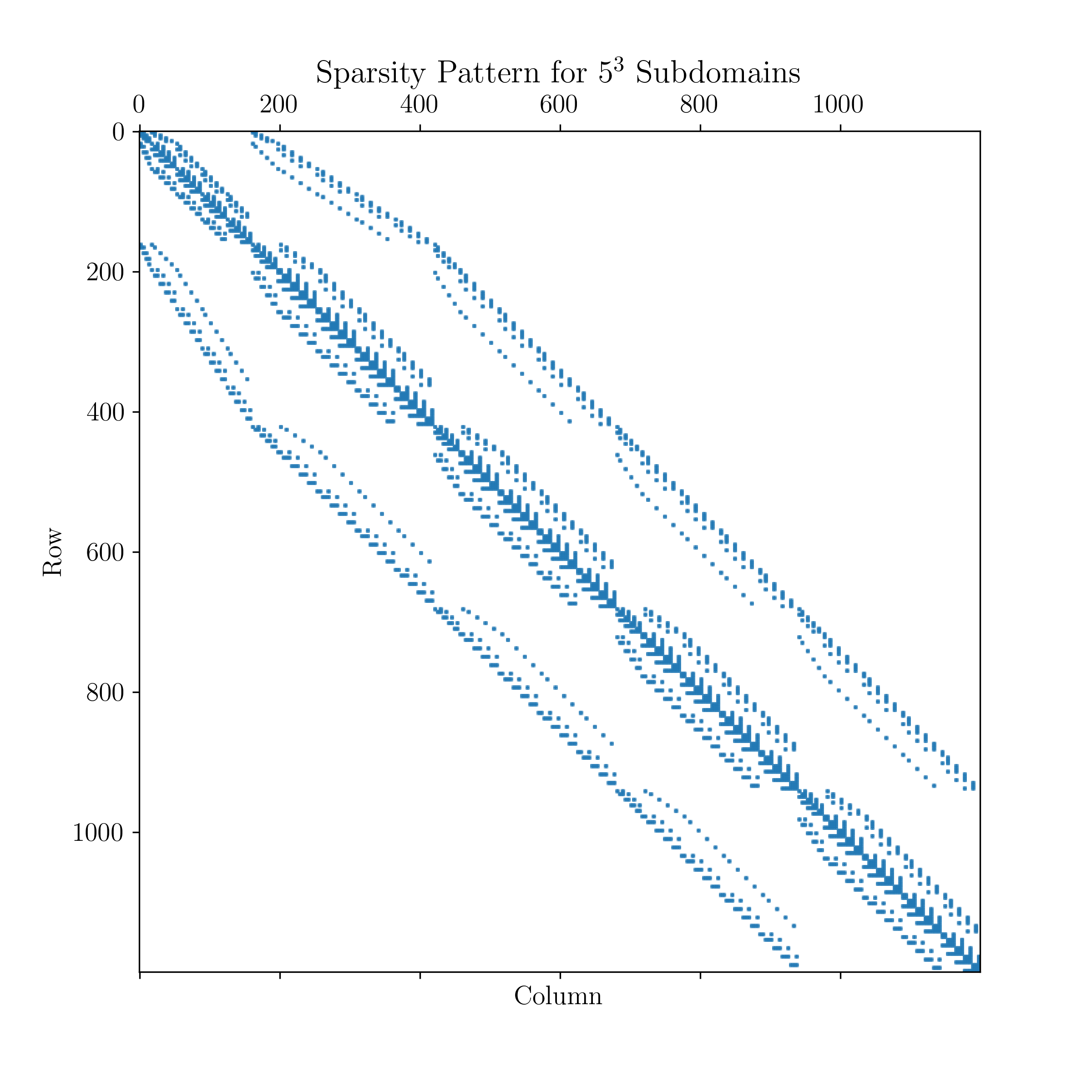}
\caption{The sparsity pattern of a problem discretized with $5 \times 5 \times 5$ boxes. Some faces have up to 10 faces on the adjacent boxes that are not on domain boundary $\partial \Omega$, but this is independent of the total number of discretization points $N$.}
\label{fig:sparsity5}
\end{figure}

\subsection{Asymptotic Cost in 3D}

The build stage has an asymptotic cost of
\begin{equation}
    O\left(N^2\right) \text{ (sparse factorization)} + O\left(p^6 N\right) \text{ (forming blocks of } \mtx{T} \text{).} \nonumber
    \label{eq:asymptoticbuild}
\end{equation}
Based on our numerical experiments the factorization of $\mathbf{T}$ is the most time-consuming portion of our method. Since we are using a multilevel solver for a 3D PDE, it has an asymptotic cost of $O(N^2)$ \cite[Chapter~25]{martinsson2019fast}. However, this cost is unaffected by the polynomial order $p$ of our discretization. Other numerical methods such as finite differences with large stencils have a significant prefactor cost in $p$ that can greatly increase the scaling constant in front of this estimate \cite[Chapter~20]{martinsson2019fast}. Additionally, since $\mathbf{T}$ is relatively small (it has $< \frac{6N}{p} \times \frac{6N}{p}$ total entries with $\approx 11p^2$ nonzero entries per row, with an example shown in Figure~\ref{fig:sparsity5}) it can be stored in memory and reused for multiple solves.

By comparison, increasing $p$ has a substantial effect on the runtime of the DtN computations needed for blocks of $\mtx{T}$. This is because each DtN $\mathbf{T}^{(i)}$ requires the solution of a system using a large block of the spectral differentiation matrix $\mathbf{A}^{(i)}$ that is size $\approx p^3 \times p^3$. Nevertheless, we can still assemble DtNs in batch for $p \leq 15$. Even for larger $p$ where batched operations are impractical, building the DtNs in serial is relatively fast as shown in our numerical results. Since DtN assembly is linear in $N$ and can be accelerated with GPUs, a steep cost in $p$ is more tolerable for it than for the factorization of $\mathbf{T}$.

Our asymptotic cost for the solve stage is
\begin{equation}
    O\left(N^{4/3}\right) \text{ (sparse solve) } + O\left(p^6 N\right) \text{ (box solves with } \mtx{S} \text{).} \nonumber
\end{equation}
If we stored the solution operators in memory then the box solves are only $O(p^2N)$, but in practice it is preferred to reformulate the solution operators and avoid having to store and potentially transfer additional memory.

\section{Numerical Experiments}
\label{sec:numerics}

In this section we present a series of numerical tests for our HPS solver on various 3D problems, investigating its accuracy and computational performance. Numerical experiments were run on two Intel(R) Xeon(R) Gold 6254 CPUs at 3.10GHz with 768GB RAM, and an Nvidia V100 GPU with 32GB of memory. Unless otherwise noted, the reported errors are the relative $\ell_2$-norm between the computed solution and the known solution as $\|\mtx{u} - \mtx{u}_{\rm true}\|_2 / \|\mtx{u}_{\rm true}\|_2$.

\begin{figure}[t]
Relative Errors for the Poisson and Helmholtz (10 ppw) Equations
\centering
\includegraphics[width=390pt,trim={78pt 0pt 90pt 25pt},clip]{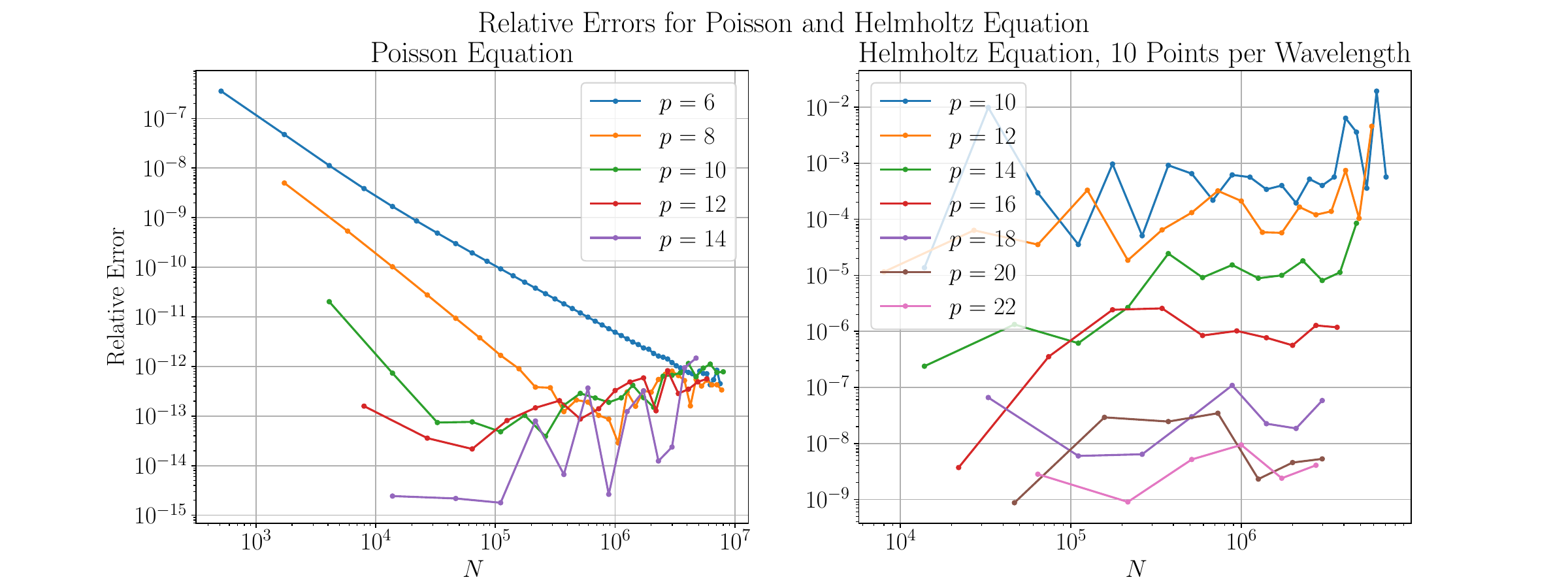}
\caption{(\textit{left}) Accuracy of our solver for the Poisson Equation. (\textit{right}) Accuracy for the Helmholtz equation fixed to 10 points-per-wavelength. Both cases feature $h$-refinement with each line marking a fixed $p$.}
\label{fig:phaccuracy}
\end{figure}

\subsection{Accuracy results - Poisson and Helmholtz Equations} \label{ssec:phaccuracy}

We investigate the accuracy of our HPS solver on the homogeneous Poisson equation,
\begin{eqnarray}
    \Delta u &=& 0 \text{   on } \Omega = [-1.1, 0.1] \times [1, 2] \times [1.2, 2.2] \nonumber \\
    u(\mathbf{x}) &=& \frac{1}{4 \pi \|\mathbf{x}\|_2} \text{   on } \partial \Omega.
    \label{eqn:poisson}
\end{eqnarray}
Here the Dirichlet boundary condition is defined so that it has a solution in its Green's function. We compute numerical solutions to this equation across a range of block sizes $h$ and polynomial orders $p$, comparing their relative $\ell^2$-norm error to the Green's function. The accuracy results are in Figure \ref{fig:phaccuracy}. We notice convergence rates of $\approx h^7$ for $p=6$ and $\approx h^{12}$ for $p=8$ until the relative error reaches $\approx 10^{-13}$, at which point it flattens (for larger $p$ the error reaches $10^{-13}$ before a convergence rate can be fairly estimated). Subsequent refinement does not consistently improve accuracy, but it does not lead to an increase in error beyond $\approx 10^{-12}$ either. This suggests our solver is relatively stable.

\begin{figure}[t]
\centering
Relative Errors for Helmholtz Equation with Fixed $\kappa$
\includegraphics[width=390pt,trim={78pt 0pt 90pt 22pt},clip]{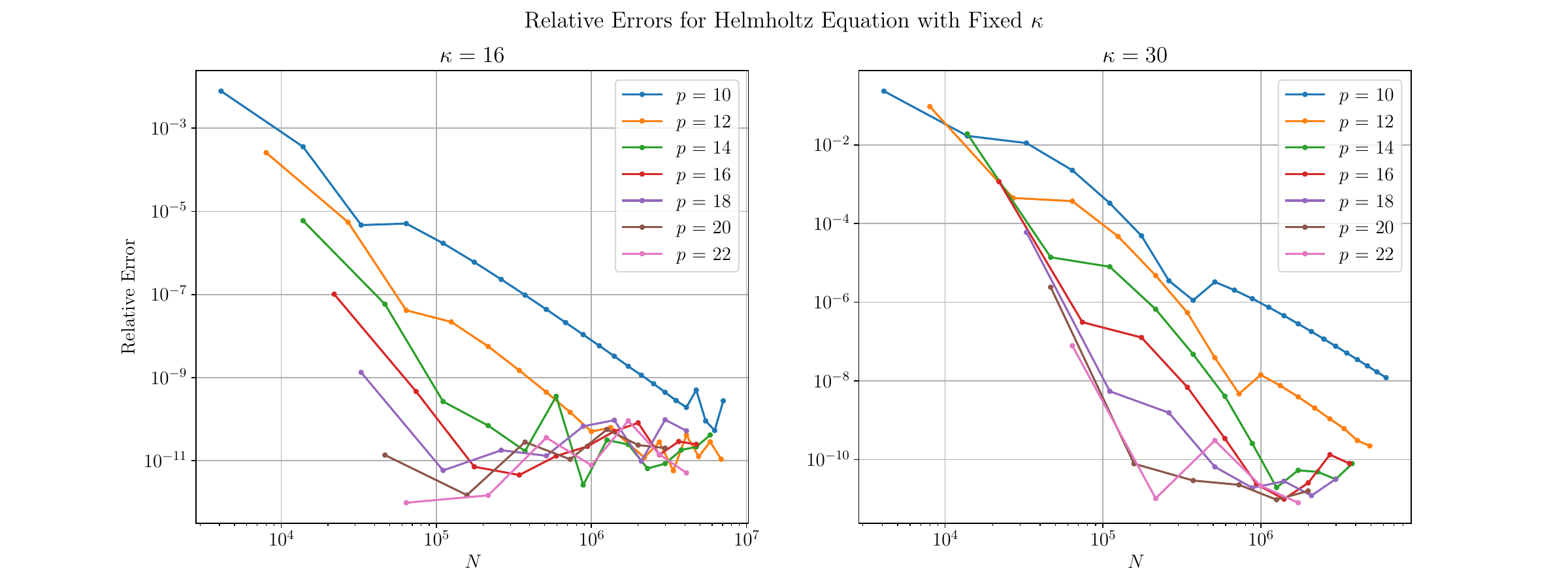}
\caption{Accuracy of our solver for the Helmholtz equation set to wavenumbers $\kappa=16$ ($\approx 2.5$ wavelengths on the domain) and $\kappa=30$ ($\approx 4.5$ wavelengths).}
\label{fig:haccuracy}
\end{figure}

To assess our method's accuracy on oscillatory problems we also evaluate it on the homogeneous Helmholtz equation with its Green's function used as a manufactured solution with wavenumber $\kappa$,
\begin{eqnarray}
    \Delta u  + \kappa^{2} u &=& 0 \text{   on } \Omega = [-1.1, 0.1] \times [1, 2] \times [1.2, 2.2] \nonumber \\
    u(\mathbf{x}) &=& \frac{\cos (\kappa |\mathbf{x}|)}{4 \pi |\mathbf{x}|} \text{   on } \partial \Omega.
    \label{eqn:helmholtz}
\end{eqnarray}
In Figure (\ref{fig:haccuracy}) we show relative $\ell^2$-norm errors for our solver on Equation (\ref{eqn:helmholtz}) with fixed wavenumbers $\kappa = 16$ and $\kappa=30$, which correspond to roughly 2.5 and 4.5 wavelengths on the domain, respectively. We see comparable convergence results to the Poisson equation, though the lower bound on relative error is somewhat higher for each problem ($\approx 10^{-11}$ for $\kappa=16$ and $\approx 10^{-10}$ for $\kappa=30$). However, increasing the resolution further does not cause the error to diverge. Along with the results for Equation (\ref{eqn:poisson}), this suggests our method remains stable when problems are over-resolved. Figure \ref{fig:phaccuracy} shows the relative error for the Helmholtz equation with a varying $\kappa$, set to provide $\approx 10$ discretization points per wavelength. We see higher $p$ corresponds to a better numerical accuracy, while increasing $h$ does not improve results (since the wavenumber is also increased), but errors remain relatively stable. Since our method maintains accuracy while increasing the total degrees of freedom linearly with $\kappa^3$, it does not appear to suffer from the pollution effect.

\begin{figure}
\centering
Gravity Helmholtz
\begin{subfigure}{.5\textwidth}
  \centering
  {\footnotesize $\kappa = 20, p = 12, h = 1/4$}
  \includegraphics[width=\textwidth,trim={5cm 7cm 4cm 8cm},clip]{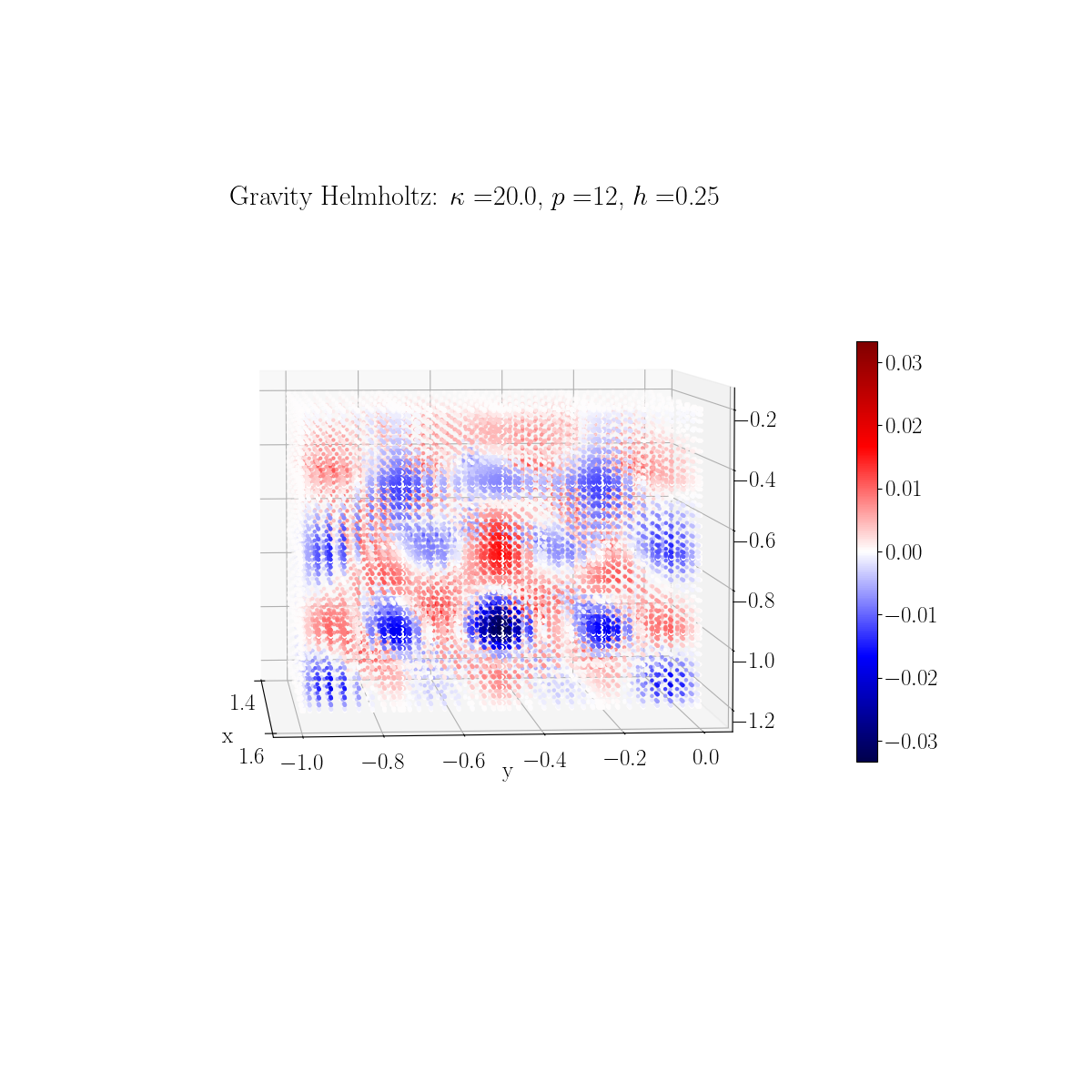}
  \caption{$\approx 64000$ discretization points.}
  \label{fig:gravityp9}
\end{subfigure}%
\begin{subfigure}{.5\textwidth}
  \centering
  {\footnotesize $\kappa = 20, p = 21, h = 1/3$}
  \includegraphics[width=\textwidth,trim={5cm 7cm 4cm 8cm},clip]{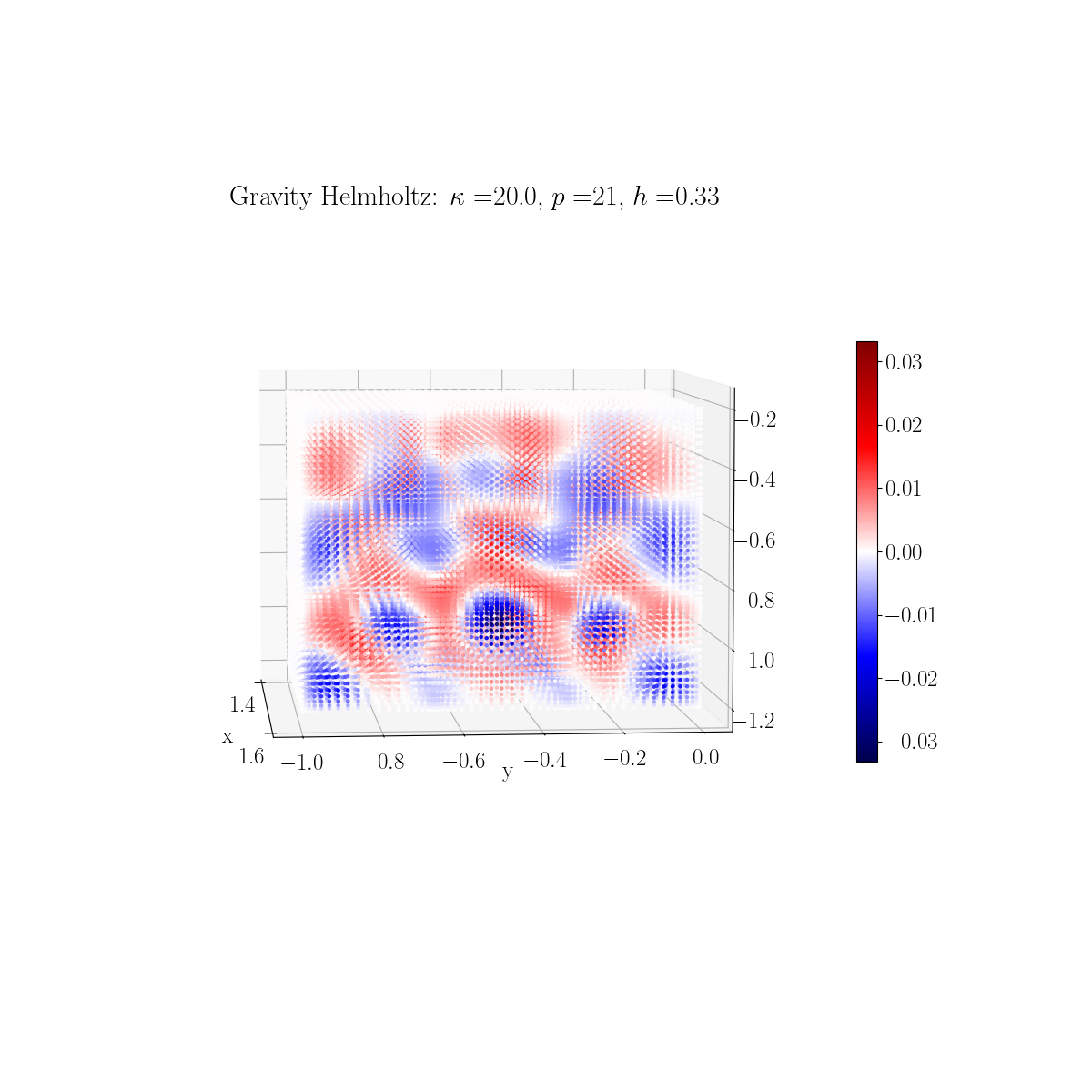}
  \caption{$\approx 185000$ discretization points.}
  \label{fig:gravityp21}
\end{subfigure}
\caption{Plots of solutions to the gravity Helmholtz equation as described in Equation~(\ref{eqn:gravityhelmholtz}). A cross-section has been taken through the $x$-axis. We see consistent results across different $p$ and $h$.}
\label{fig:gravityscatter}
\end{figure}

\begin{figure}[t]
\centering
$p$-Refinement of Gravity Helmholtz Equation
\includegraphics[width=300pt,trim={0pt 0pt 0pt 50pt},clip]{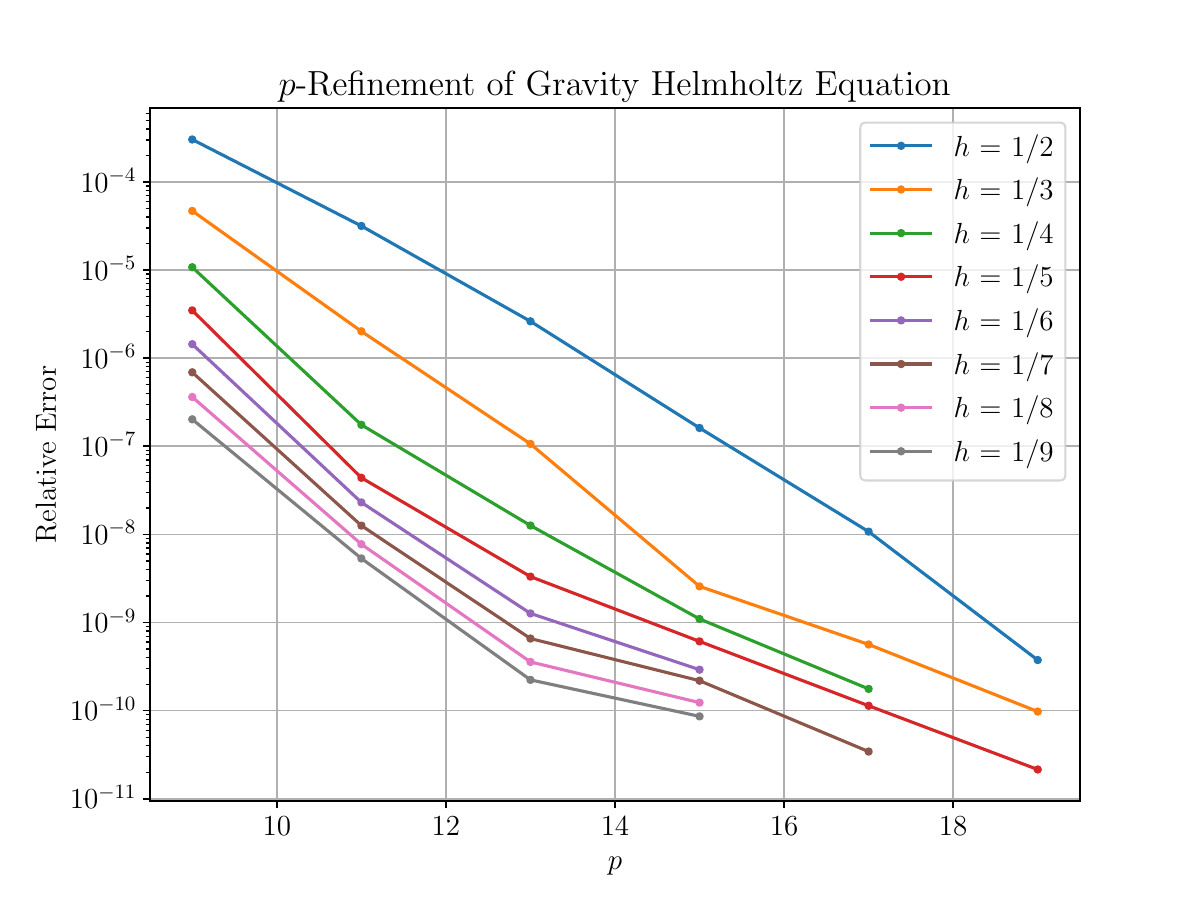}
\caption{Convergence of our solver on the gravity Helmholtz equation with $\kappa=12$, relative to a highly-resolved solution. Here $h$ is fixed and resolution is increased through $p$-refinement.}
\label{fig:gravityN}
\end{figure}

We then consider the gravity Helmholtz equation in an example similar to one in \cite{fortunato2021ultraspherical},
\begin{eqnarray}
    \Delta u  + \kappa^2 (1-x_3) u &=& -1 \text{   on } \Omega = [1.1, 2.1] \times [-1, 0] \times [-1.2, -0.2] \nonumber \\
    u(\mathbf{x}) &=& 0 \text{   on } \partial \Omega.
    \label{eqn:gravityhelmholtz}
\end{eqnarray}
This equation adds a spatially-varying component to the wavenumber variation based on $x_3$ that is analogous to the effect of a gravitational field, hence the name \cite{barnett2015high}. Since $1 < 1-x_3 < 2.2$, the variation is greater through $\Omega$ than in the constant Helmholtz equation with the same $\kappa$. Unlike the previous Poisson and Helmholtz equations, equation (\ref{eqn:gravityhelmholtz}) does not have a manufactured solution. It also produces relatively sharp gradients near the domain boundaries since there is a zero Dirichlet boundary condition, and it has a variable coefficient in spatial coordinates (specifically $x_3$). We investigate our solver's accuracy by computing the relative $\ell^2$-norm errors of our results to an over-resolved solution at select points. Figure \ref{fig:gravityscatter} shows the solution of this problem using select $p$ and $h$ values, while Figure \ref{fig:gravityN} shows these errors in the case of $p$-refinement for a range of values of $h$. We can see our solver shows convergence for each $h$.

\subsection{Speed and scaling}
\label{ssec:speedscaling}

\begin{figure}[t]
\centering
Build Stage Times
\includegraphics[width=390pt,trim={78pt 0pt 90pt 5pt},clip]{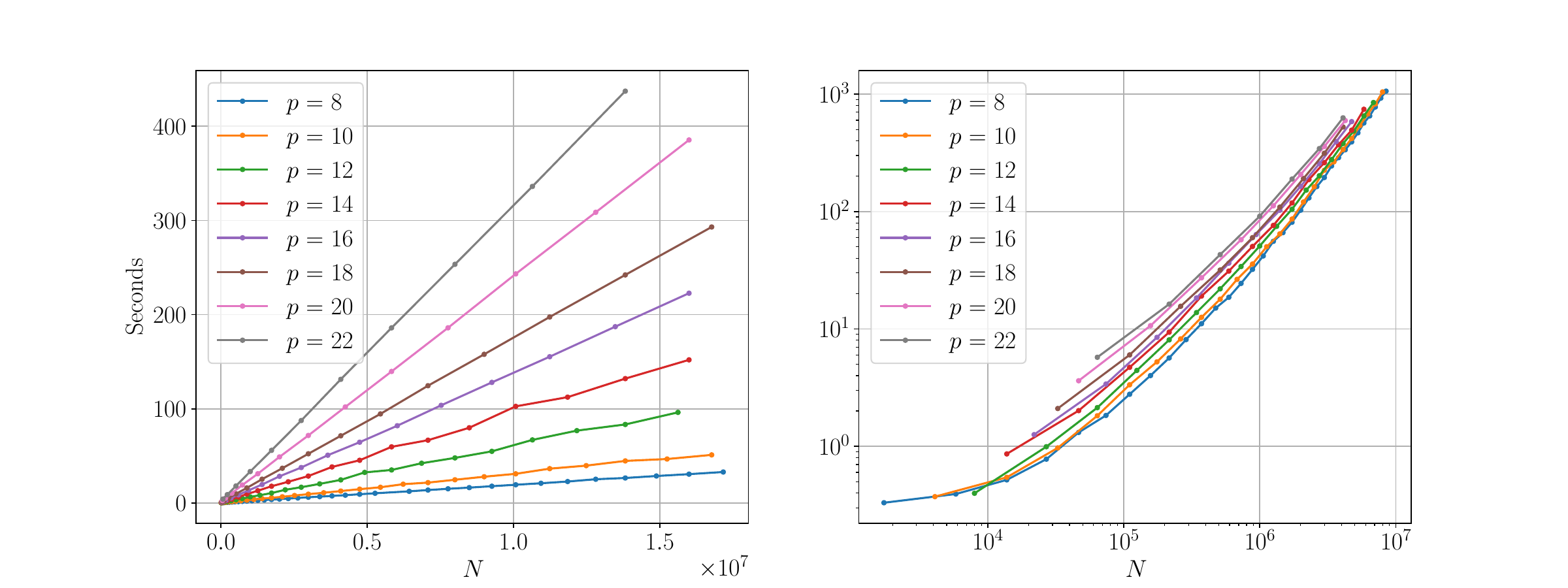}
\caption{Time to formulate the DtN maps necessary for the sparse system $\mtx{T}$ (\textit{left}) and factorize $\mtx{T}$ (\textit{right}).}
\label{fig:buildTime}
\end{figure}

To evaluate the computational complexity of our HPS solver, we measure the execution times of four components: the assembly of the DtN operators used for the build stage, the factorization time for our sparse matrix constructed by the DtNs, the solve time for the factorized system, and the time for batched solves of the box subdomain interiors. We know the asymptotic cost of constructing the DtN operators is $O \left(p^6 N \right)$ as shown in Equation (\ref{eq:asymptoticbuild}). Figure (\ref{fig:buildTime}) illustrates the run times for the build stage in solving the Poisson and Helmholtz equations detailed in Section \ref{ssec:phaccuracy}, using batched linear algebra routines on an Nvidia V100 GPU. As expected, we observe a linear increase in runtime with $h$-refinement (i.e. a higher $N$ with fixed $p$) for the DtN assembly. The constant for the linear scaling in $h$, though, increases sharply with higher values of $p$.

\begin{figure}[t]
\centering
$\mtx{T}$ Solve Times
\includegraphics[width=390pt,trim={78pt 0pt 90pt 5pt},clip]{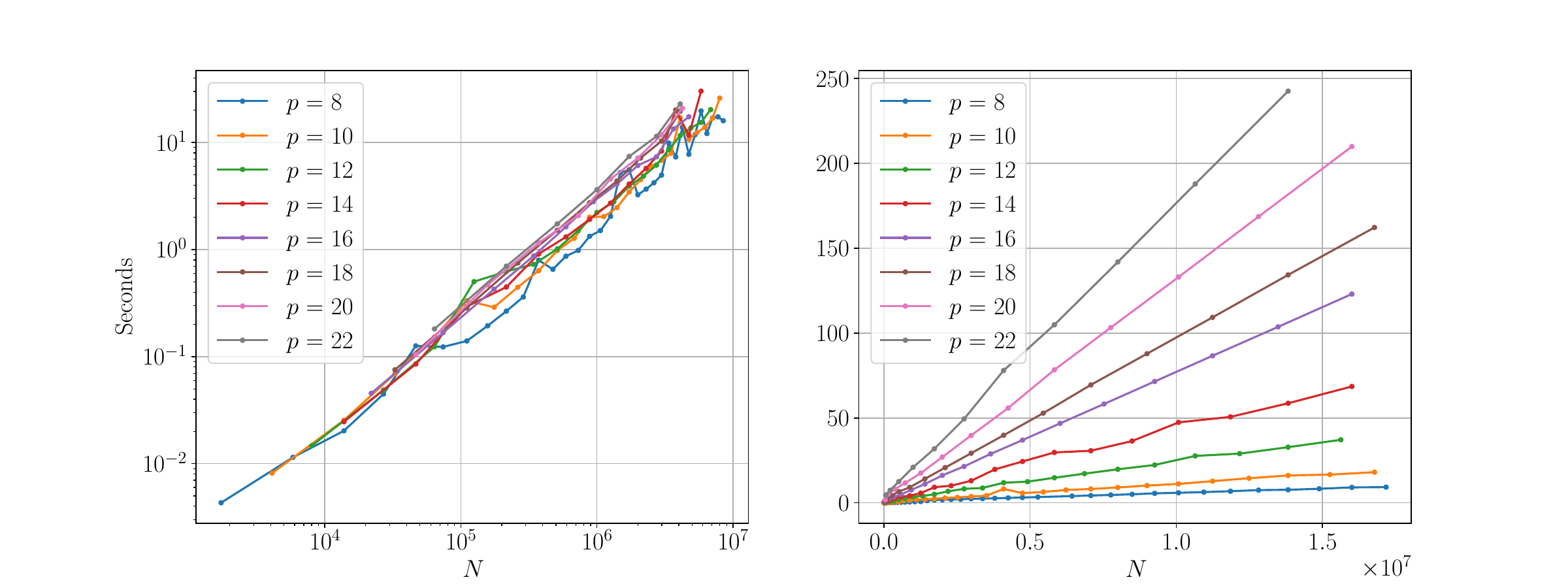}
\caption{Time to solve a problem using the factorized sparse system $\mtx{T}$ (\textit{left}) and apply the interior box solves (\textit{right}).}
\label{fig:solvetime}
\end{figure}

In Figure \ref{fig:buildTime} we also see the time to factorize the resulting sparse matrix for the build stage, using MUMPS. The expected asymptotic cost is $O(N^2)$. However, our results closely resemble a lower cost of $O(N^{3/2})$ for $N$ up to $\approx$ 8 million. Overall the factorization time dominates the build stage of our method, given its more costly scaling in $N$ and the reliance on CPUs. Lastly in Figure \ref{fig:solvetime} we see the runtime for the solves of factorized $\mtx{T}$ for the box boundaries and the batched box interior solves, respectively. Although the box boundary solve has a higher asymptotic cost in $N$, in practice it has a shorter runtime, especially for large $p$. The box interior solves follow a linear scaling similar to the batched DtN construction.

\subsection{Curved and non-rectangular domains}
\label{ssec:curved}

\begin{figure}
\centering
Helmholtz Equation on Non-Rectangular Domains
\begin{subfigure}{.5\textwidth}
  \centering
  {\footnotesize $\kappa = 26, p = 18, h = 1/4$}
  \includegraphics[width=\textwidth,trim={4cm 5cm 4cm 8cm},clip]{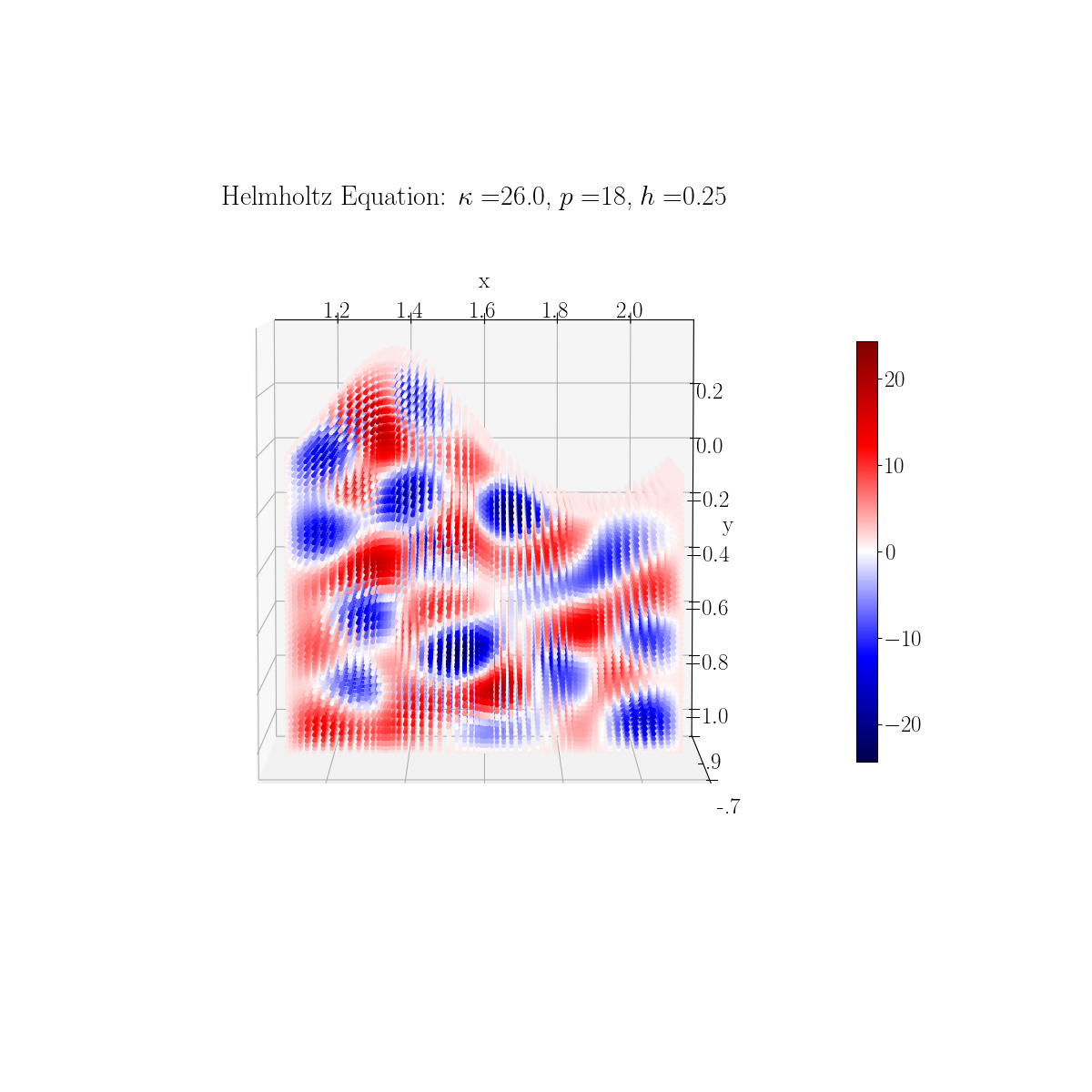}
  \caption{On a curved domain with a wavy edge, $\approx 2.6 \times 10^5$ discretization points.}
  \label{fig:curved}
\end{subfigure}%
\begin{subfigure}{.5\textwidth}
  \centering
  {\footnotesize$\kappa = 8, p = 14, h = 1/16 \times 1/2 \times 1/2$}
  \includegraphics[width=\textwidth,trim={4cm 5cm 4cm 8cm},clip]{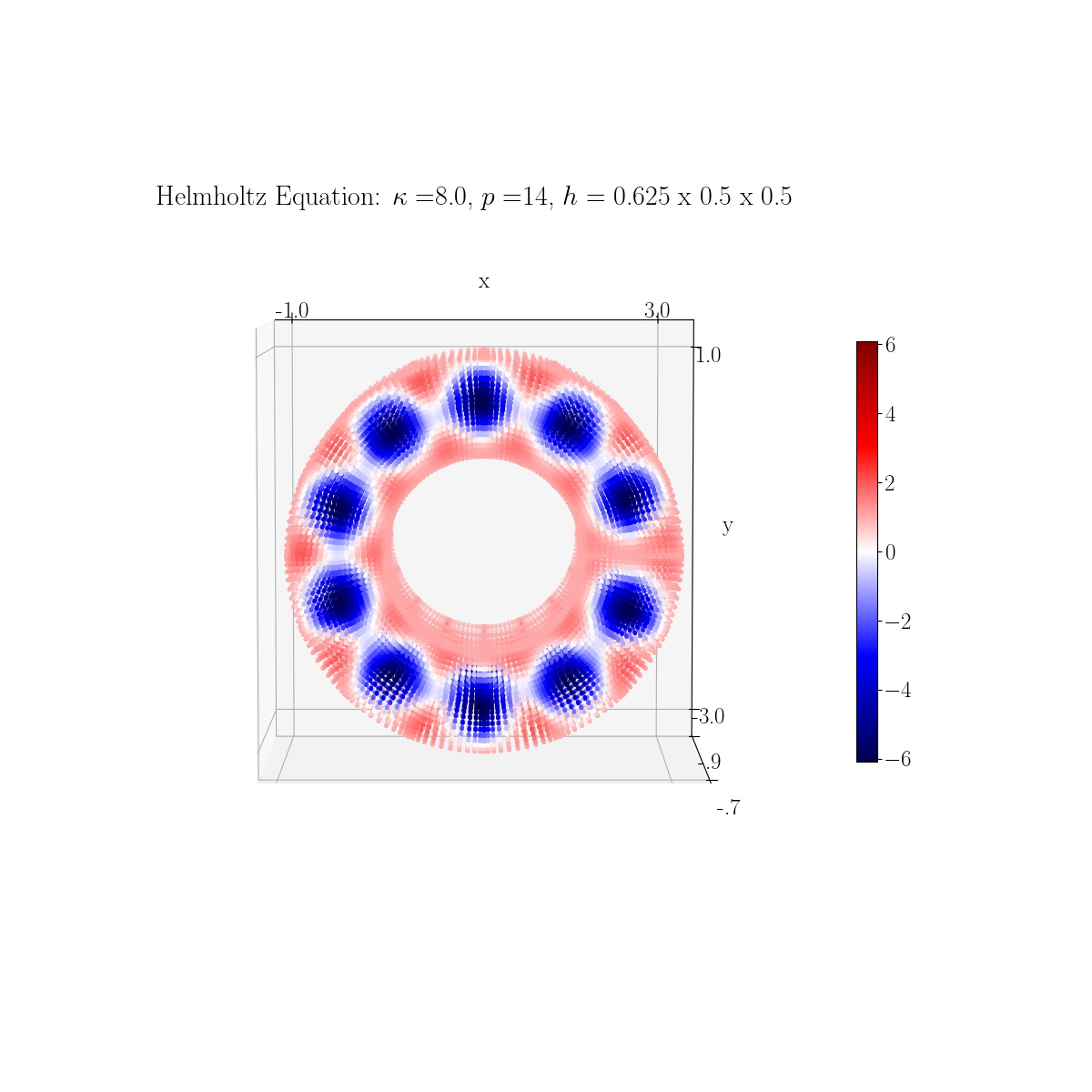}
  \caption{On an annulus, with $\approx 6.9 \times 10^4$ discretization points.}
  \label{fig:annulus}
\end{subfigure}
\caption{Helmholtz equation on two non-rectangular domains, both with a Dirichlet boundary condition of 1. A cross section has been taken along the $z$ axis.}
\label{fig:curvedDomain}
\end{figure}

\begin{figure}[t]
\centering
Relative Errors for Helmholtz Equation on a Curved Domain
\includegraphics[width=390pt,trim={78pt 0pt 90pt 22pt},clip]{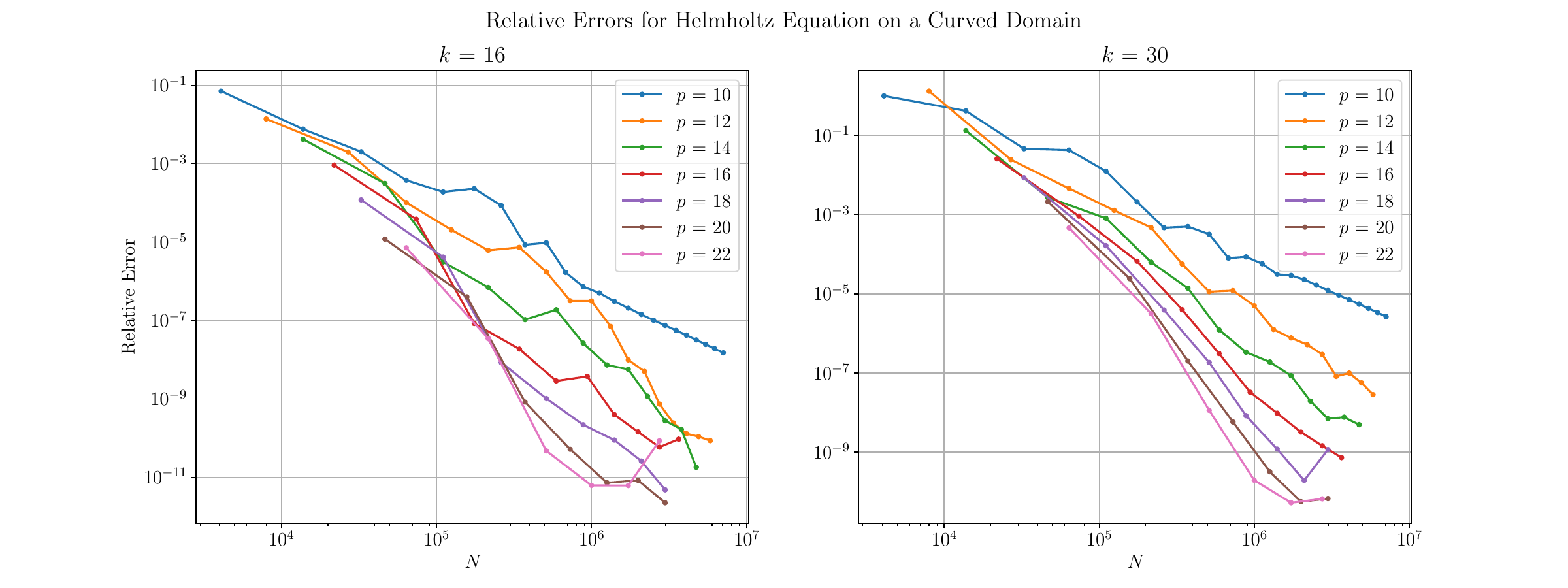}
\caption{Accuracy of our solver for the Helmholtz equation on a curved domain shown in Fig. \ref{fig:curved} set to wavenumbers $\kappa=16$ ($\approx 2.5$ wavelengths on domain) and $\kappa=30$ ($\approx 4.5$ wavelengths).}
\label{fig:curvedaccuracy}
\end{figure}

The HPS method has been extended to other problems such as surface PDEs \cite{fortunato2024high}, inverse scattering problems \cite{borges2017high}, and more generally non-rectangular domains. We can apply our HPS solver to non-rectangular domain geometries through the use of an analytic parameterization between the domain we wish to model, such as a sinusoidal curve along the $y$-axis shown in Figure~\ref{fig:curved}, and a rectangular reference domain. These parameter maps extend the versatility of the solver, but feature multiple challenges in their implementation: they require variable coefficients for our differential operators; they may use mixed second order differential terms, i.e. $\partial^2 u / (\partial x_i \partial x_j)$ where $i \neq j$; and (depending on the domain) may have areas that are near-singular, such as around the sharp point $(x, y) \approx (2, 0)$ in Figure~\ref{fig:curved}. We investigate the accuracy of our solver for the Helmholtz equation as shown in Equation~(\ref{eqn:helmholtz}), but now on the sinusoidal domain
\begin{equation}
    \Psi = \left\{ x_1, \frac{x_2}{\psi(x_1)}, x_3 \text{ for } (x_1, x_2, x_3) \in \Omega \right\} \text{ where } \psi(z) = 1 - \frac{1}{4}\sin(6z).
    \label{eqn:curveddomain}
\end{equation}
We map Eq.~(\ref{eqn:helmholtz}) on $\Psi$ to a cubic reference domain $\Omega = [1.1,2.1] \times [-1, 0] \times [-1.2, 0.2]$ using the parameter map
\begin{eqnarray}
    -\frac{\partial^2 u}{\partial x_1^2} - \left(\left(\frac{\psi'(x_1)x_2}{\psi(x_1)}\right)^2 + \psi(x_1)^2 \right)\frac{\partial^2 u}{\partial x_2^2} \nonumber - \frac{\partial^2 u}{\partial x_3^2} - 2\frac{\psi'(x_1)x_2}{\psi(x_i)} \frac{\partial^2 u}{\partial x_1 \partial x_2} \\
    - \frac{\psi''(x_1)x_2}{\psi(x_i)} \frac{\partial u}{\partial x_2} - \kappa u = 0, \text{ } (x_1, x_2, x_3) \in \Omega.
    \label{eq:parametermap}
\end{eqnarray}
We test this problem with a manufactured solution in the Green's function of the Helmholtz equation, with corresponding Dirichlet data matching our test in Section~\ref{ssec:phaccuracy}. Convergence studies for this problem with $\kappa =16$ and $\kappa =30$ are shown in Figure \ref{fig:curvedaccuracy}. We observe a similar convergence pattern to the Helmholtz test on a cubic domain, although the speed of convergence is somewhat slower. A plot of this problem for wavenumber $\kappa = 26$ with a zero Dirichlet boundary condition and body load $f=1$ is shown in Figure \ref{fig:curved}, to provide an example that does not have a manufactured solution.

We also investigate our solver's accuracy on the Helmholtz equation when the domain is a three-dimensional annulus, with an example with zero Dirichlet boundary condition and $f = 1$ shown in Figure \ref{fig:annulus}. We use an analytic parametrization to map this problem to a rectangular reference domain with higher resolution in the $y$-axis and a periodic boundary condition. Relative errors in $p$ and $N$ are comparable to those for the problem on the sinusoidal domain.

\subsection{Parabolic PDEs with time stepping}

Consider a convection-diffusion equation modeling the concentration of a contaminant $u$ in a fluid, over time interval $[0,T=5]$ on $\Omega = [-0.5, 0.5]^3$, with diffusivity $\kappa = 10^{-4}$ and circular horizontal velocity field $\mathbf{b}(\mathbf{x})$:
\begin{align}
    \label{eqn:convdiff}
    \mathbf{b}(x_1,x_2,x_3) &= \left(-\cos (x_1) \sin(x_2) x_3, \sin(x_1)\cos(x_2)x_3, 0\right) \nonumber \\
    \frac{\partial u(\mathbf{x},t)}{\partial t} &= \kappa \Delta u - \nabla \cdot (\mathbf{b}(\mathbf{x}) u) = Au, &(x,t) \in \Omega \times [0, T], \nonumber \\
    u(\mathbf{x}, t) &= 0, &(x, t) \in \partial \Omega \times [0,T].
\end{align}
Here $A$ is the linear partial differential operator where $Au = (\kappa \Delta u - \nabla \cdot (\mathbf{b} u))$. Equation~(\ref{eqn:convdiff}) is a parabolic PDE that requires both temporal and spatial discretization. We will apply a temporal discretization through the Crank-Nicolson method, which gives us
\begin{align}
    & \frac{u^{n+1} - u^n}{\Delta t} = \frac{1}{2} A u^{n+1} + \frac{1}{2} A u^n \nonumber \\
    \implies& \left(I - \frac{\Delta t}{2}A\right)u^{n+1} = \left(I + \frac{\Delta t}{2}A\right)u^n.
    \label{eqn:backwardEuler}
\end{align}
The left-hand-side of (\ref{eqn:backwardEuler}) is itself a linear elliptic partial differential operator, thus we may discretize it using our HPS method and produce both box solution and DtN maps as well as the system $\mtx{T}$. Moreover, the same PDO is used for every time step in the computation. Thus the build stage of our method need only be applied once - we can use the same factorized $\mathbf{T}$ for every time step, supplying $\mathbf{u}^n$ from the previous timestep as the body load. However, we must also multiply $\mathbf{u}^n$ by a similar but distinct operator, $I + (\Delta t / 2)A$. This is applied to both the box solves (by setting $\mathbf{f} = (\mathbf{I} + (\Delta t / 2)\mathbf{A})[I,:] \mathbf{u}^n$ for our interior box solve), and the right-hand-side for the sparse system solve in $\mtx{T}$. Previous works such as \cite{chen2024fast, babb2020hps} have explored this concept for using HPS to solve time-dependent problems.

\begin{figure}
\centering
\begin{subfigure}{.33\textwidth}
  \centering
  \includegraphics[width=\textwidth,trim={1.5cm 1.5cm 2cm 2cm},clip]{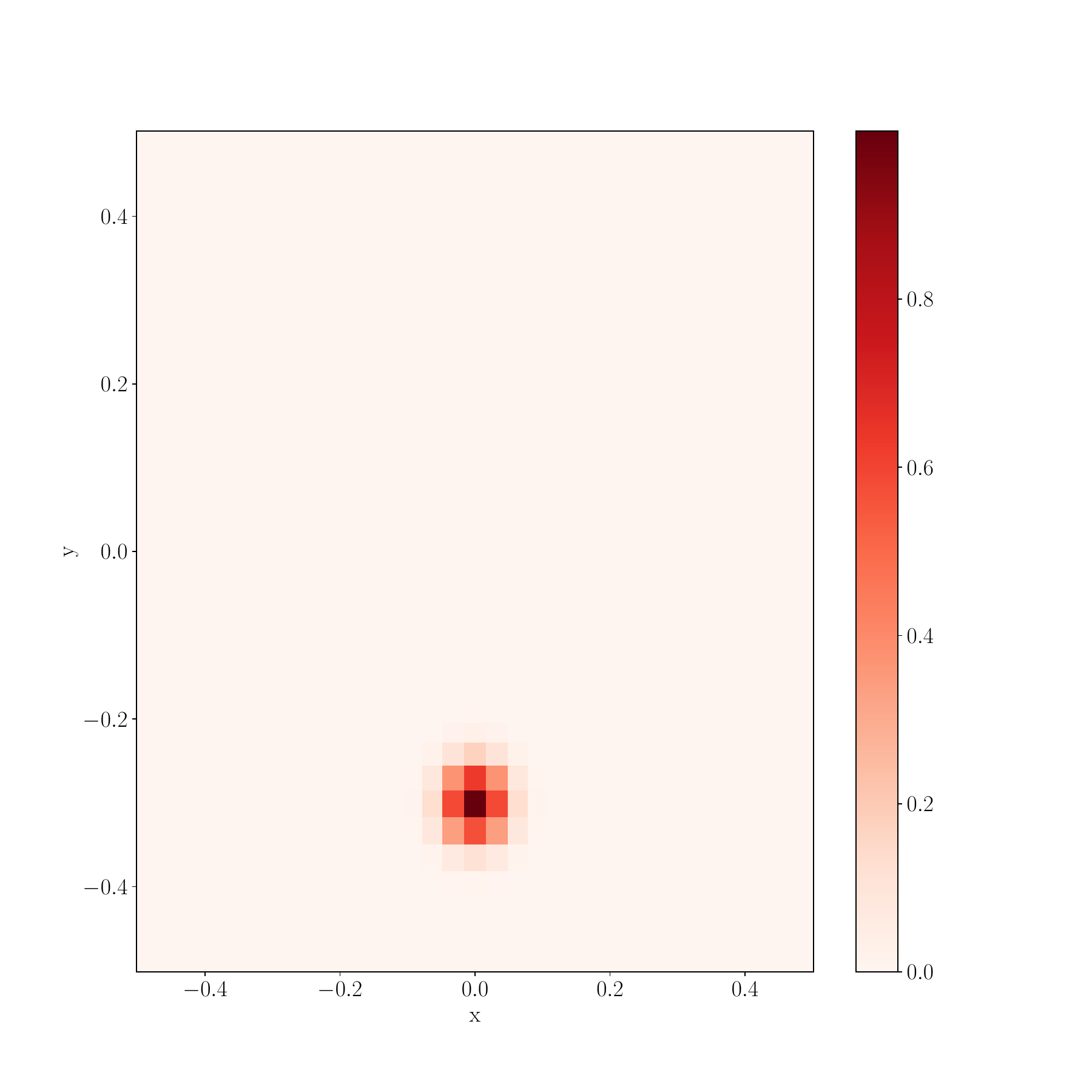}
  \caption{Initial Condition.}
\end{subfigure}%
\begin{subfigure}{.33\textwidth}
  \centering
  \includegraphics[width=\textwidth,trim={1.5cm 1.5cm 2cm 2cm},clip]{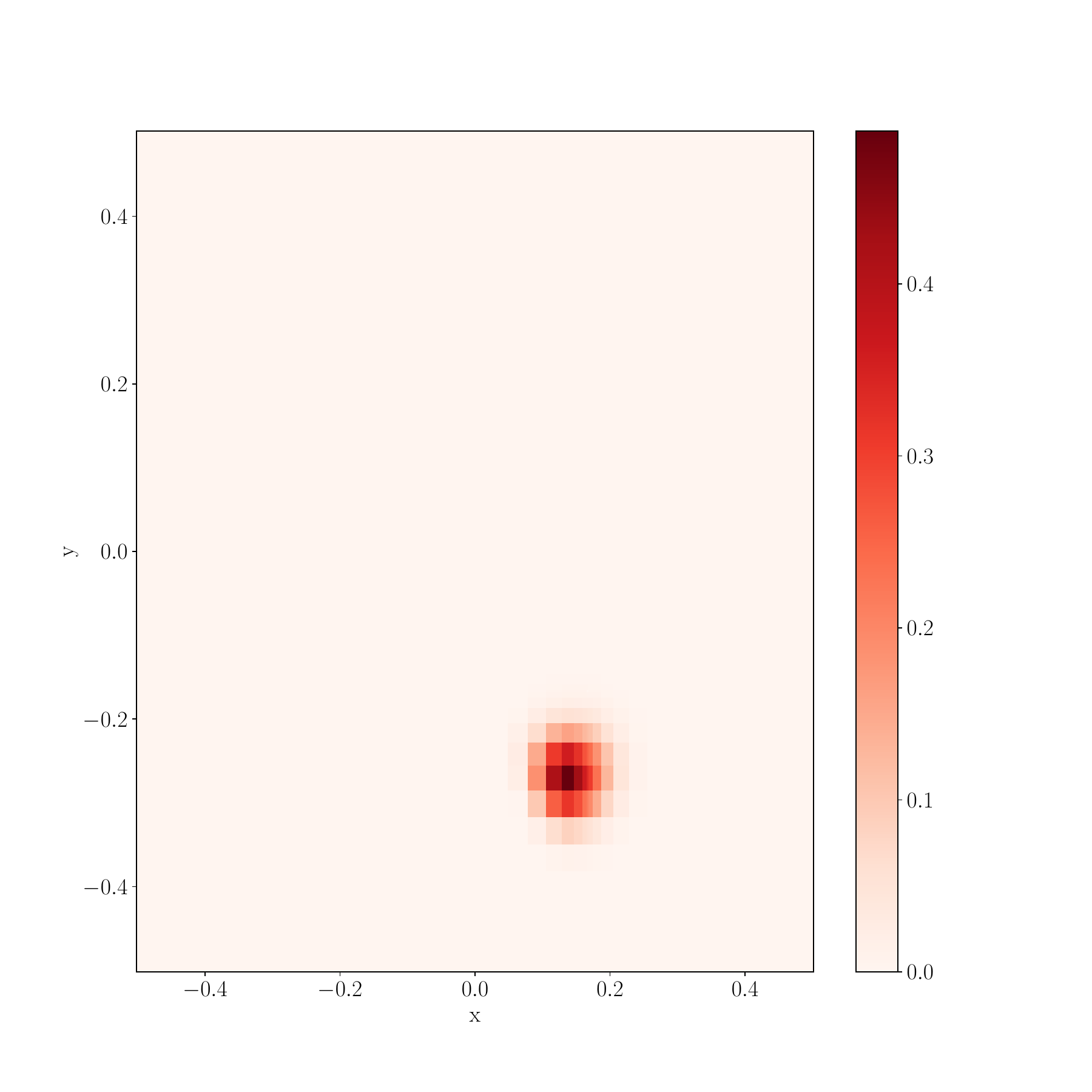}
  \caption{After 10 time steps.}
\end{subfigure}
\begin{subfigure}{.33\textwidth}
  \centering
  \includegraphics[width=\textwidth,trim={1.5cm 1.5cm 2cm 2cm},clip]{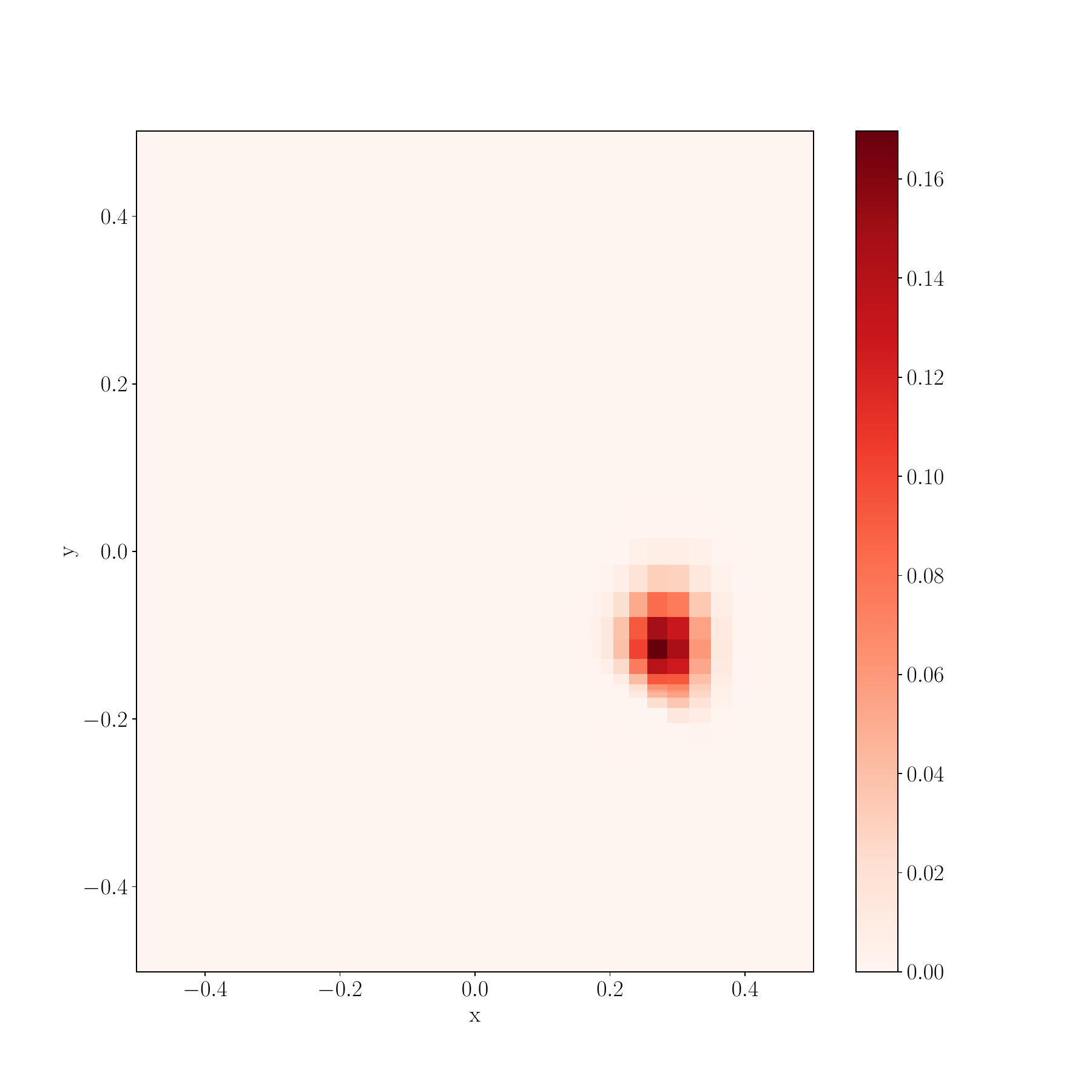}
  \caption{After 25 time steps.}
\end{subfigure}
\caption{Progression of contaminant $\mathbf{u}$ in a circular flow up to $t=2.5$. This is a cross-section of a three-dimensional space, so vertical diffusion also affects the result. We see $\mathbf{u}$ follows an expected counterclockwise pattern, continued in Figure~\ref{fig:timestepComparison}.}
\label{fig:timestepProgress}
\end{figure}

We investigate the use of our HPS method in solving the PDE in Equation~(\ref{eqn:convdiff}) with the initial condition $u^0 = \exp((x_1^2 + (x_2+0.3)^2 + x_3^2)/0.002)$. This creates a contaminant density positioned slightly below the center of $\Omega$ along the $y$ axis. Thus the contaminant $u$ will be carried by the velocity field in a counterclockwise pattern. The progress of this flow can be seen in Figure \ref{fig:timestepProgress}. As expected, the contaminant is beginning to move in the pattern determined by $\mathbf{b}(\mathbf{x})$ as the convective term.

\begin{figure}
\centering
\begin{subfigure}{.33\textwidth}
  \centering
  \includegraphics[width=\textwidth,trim={1.5cm 1.5cm 2cm 2cm},clip]{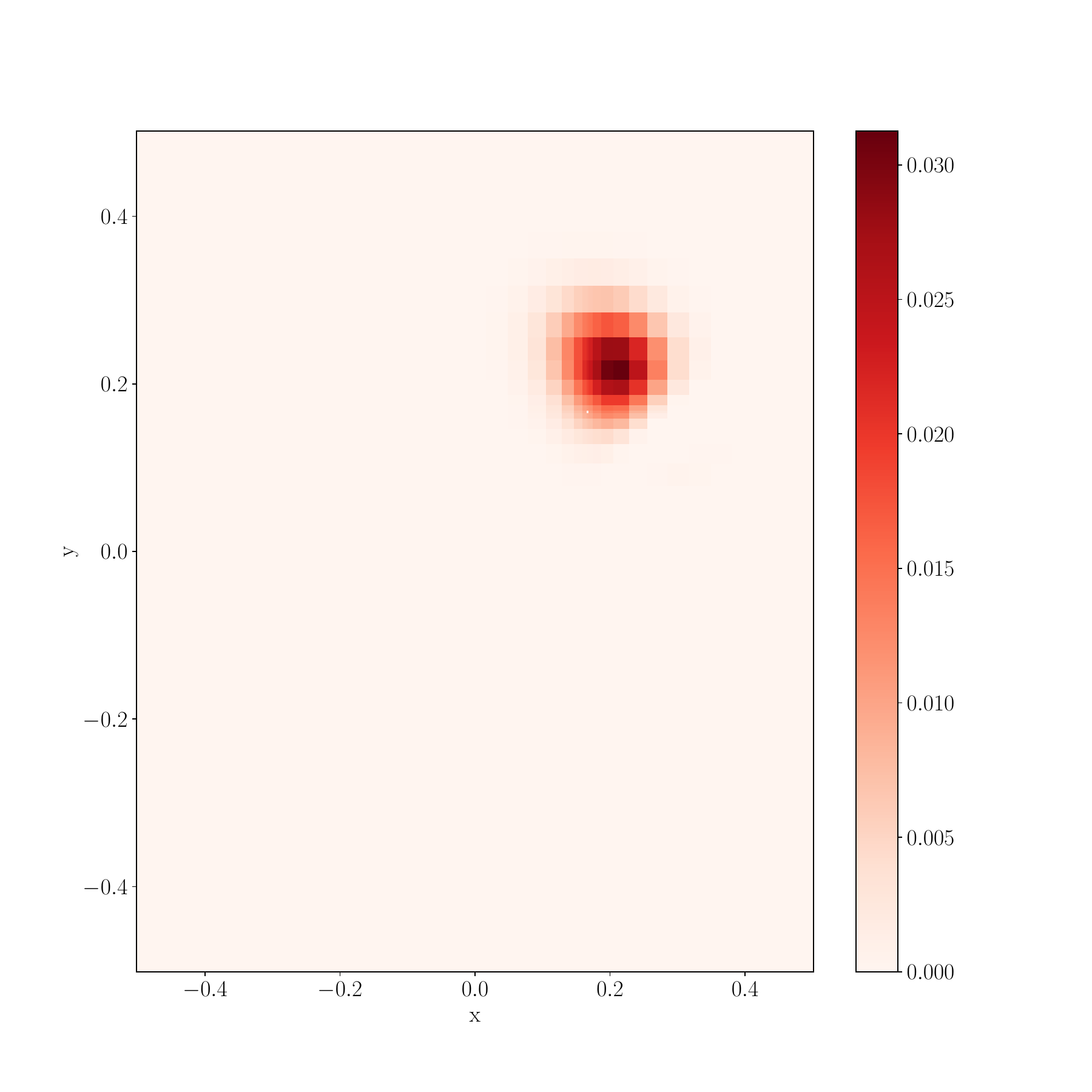}
  \caption{Higher spatial resolution.}
\end{subfigure}%
\begin{subfigure}{.33\textwidth}
  \centering
  \includegraphics[width=\textwidth,trim={1.5cm 1.5cm 2cm 2cm},clip]{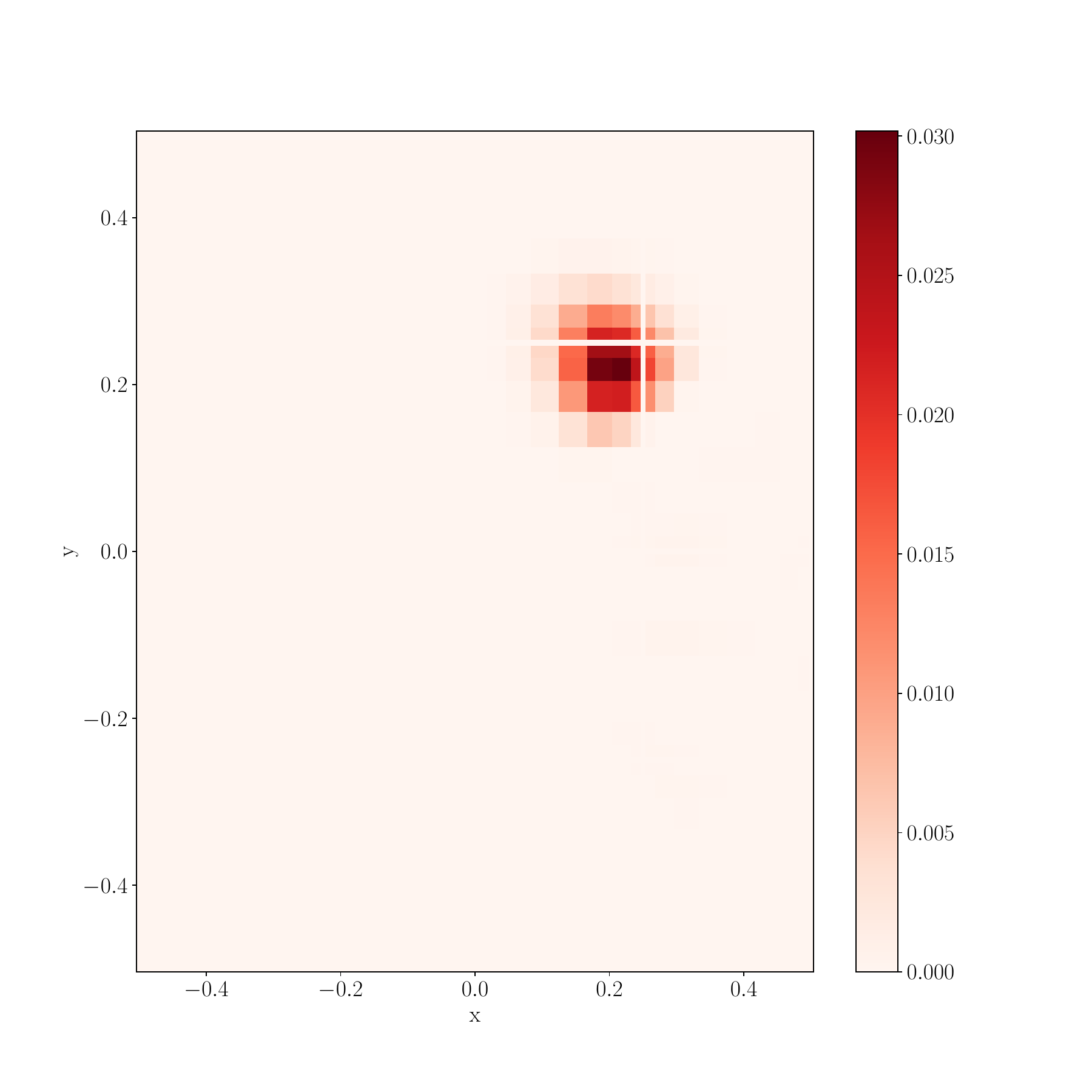}
  \caption{Low space resolution, larger timestep.}
\end{subfigure}
\begin{subfigure}{.33\textwidth}
  \centering
  \includegraphics[width=\textwidth,trim={1.5cm 1.5cm 2cm 2cm},clip]{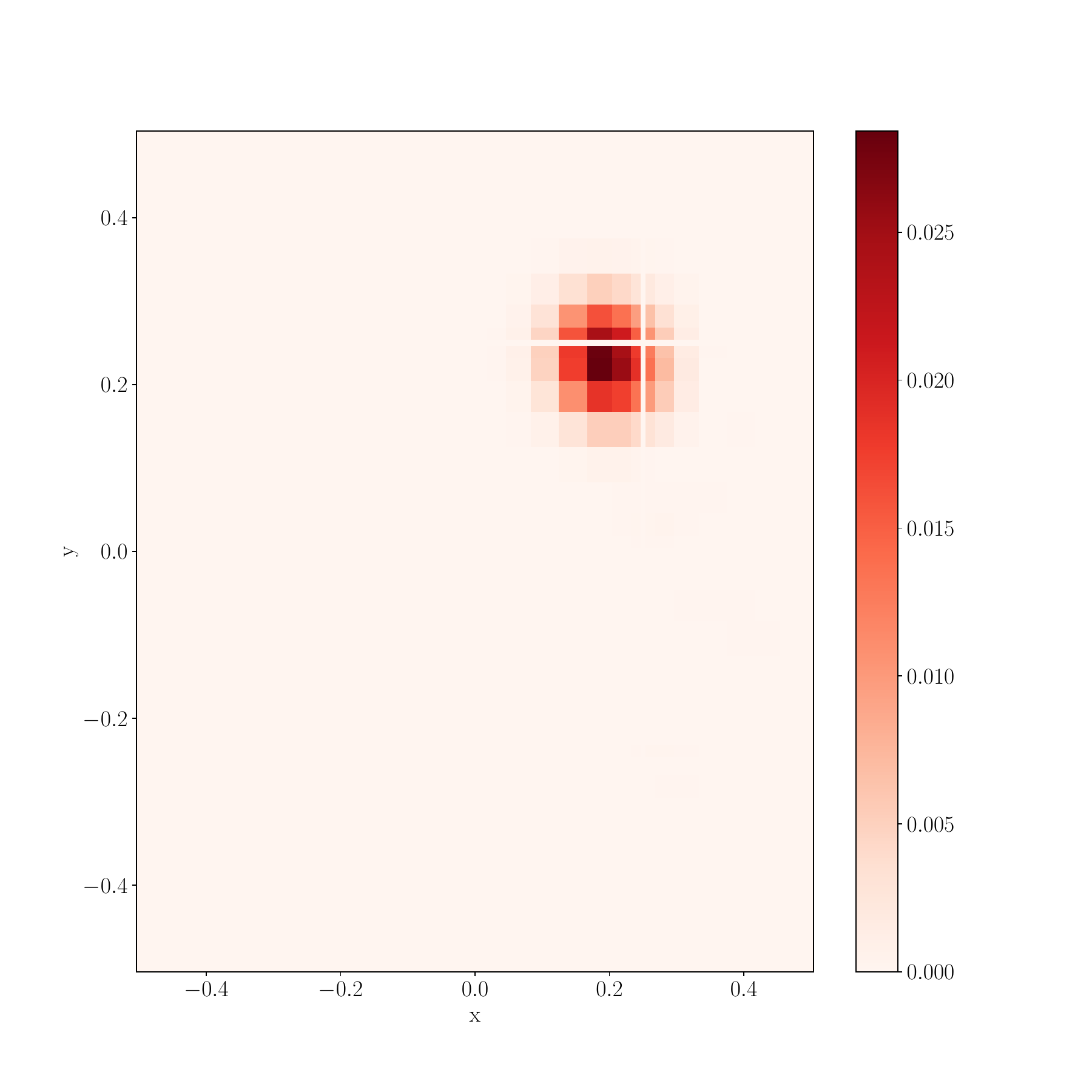}
  \caption{Smaller timestep.}
\end{subfigure}
\caption{Comparison of $\mathbf{u}$ at $t=5$ in three regimes: lower spatial resolution and larger time step in the middle, with higher space or time resolutions on either side.}
\label{fig:timestepComparison}
\end{figure}

For modeling a parabolic PDE we must consider the impact of both spatial and temporal refinement. The Crank-Nicolson method is numerically stable, but it is only a second-order method. Thus refining the size of timestep $\Delta t$ may have a more significant effect on numerical results than refining $p$ or $h$ in the HPS method. This is explored in Figure~\ref{fig:timestepComparison}. Reducing $\Delta t$ does increase the method's overall runtime, but since the same sparse system $\mathbf{T}$ is used for each step we bypass the build stage of the method, including the costly factorization of $\mathbf{T}$. In practice HPS would be better combined with a higher order time stepping scheme such as certain Runge-Kutta methods, but as long as the underlying partial differential operator is linear we may leverage the advantage of applying the build stage only once.

\section{Conclusion and Future Work}

We present a Hierarchical Poincar\'e–Steklov scheme tailored for modern hardware, combining GPU-accelerated local solves with a black-box multifrontal sparse solver to deliver both numerical stability and scalability. By casting the global interface system in a form amenable to multifrontal factorization, we decouple dense local computations from the sparse interface solve. In our two-level approach, batch-processed linear algebra on the GPU first factorizes the block-diagonal subdomain systems; then a sparse system on the remaining interfaces is assembled and factorized using black-box sparse solvers, yielding enhanced performance and robustness.

The method presented can readily be adapted to handle other versions of HPS, including those relying on impedance-to-impedance maps \cite{gillman2015spectrally} and those that use expansion coefficients to represent interface functions \cite{fortunato2021ultraspherical}. It could in principle also be used in situations where adaptivity in either the local discretization order $p$ or in the mesh size is used \cite{chipman2024fast,geldermans2019adaptive}, or with formulations that exploit internal structure in the local spectral differentiation matrices \cite{aurentz2015fast}. However, it may be a delicate matter to maintain the very high efficiency provided by batched linear algebra in our implementation.

Rank‐structured compression offers a means of reducing the cost of the global interface solve. Dense Schur complements arising from the coarse level of a sparse factorization are often rank‐compressible \cite{martinsson2019fast,2003_borm_introduction_H_matrix,2007_leborne_HLU,2013_xia_rankpatterns} and can be recovered using black-box randomized methods \cite{levitt2024linear,2011_lin_lu_ying,2011_martinsson_randomhudson,pearce2025randomizedrankstructuredmatrixcompression,2013_xia_randomized,yesypenko2025randomizedstrongrecursiveskeletonization}. 
Slab‐based domain decomposition schemes \cite{engquist2011sweeping,gander2019class,yesypenko2024slablu} are particularly compelling because they confine most operations to local subdomains, so that global coupling is restricted to a small, coarse‐scale system --- simplifying parallelism and enabling effective use of rank‐structured compression compared to nested‐dissection approaches. Pursuing these ideas will allow the HPS method to scale efficiently to larger problems while preserving robustness for complex physical systems.

\section{Code Availability}

A code repository of the HPS solver described in this paper is available here: \url{https://doi.org/10.5281/zenodo.16379968} \cite{joseph_kump_2025_16379968}.

\section{Conflicts of Interest and Acknowledgments}

The authors have no conflicts of interest to report.

This work was supported by the Office of Naval Research (N00014-18-1-2354), by the
National Science Foundation (DMS-1952735 and DMS-2313434), and by the Department of Energy
ASCR (DE-SC0022251).


\bibliographystyle{unsrt}
\bibliography{refs.bib}



\end{document}